\documentclass{article}
\usepackage{amssymb, amsthm, amsmath}
\usepackage{graphicx}
\usepackage[T2A]{fontenc}
\usepackage{float}
\usepackage[usenames]{color}
\usepackage{colortbl}
\usepackage{subfigure}
\usepackage{enumitem}
\usepackage{authblk}

\newcommand{\rp}{\mathbb R \operatorname{P}}

\newtheorem{theorem}{Theorem}
\newtheorem{lemma}{Lemma}

\newtheorem{definition}{Definition}
\newtheorem{proposition}{Statement}

\begin{document}
\sloppy
\title{Non-singular flows with twisted saddle orbit on orientable 3-manifolds}
\author{Olga Pochinka, Danila Shubin\footnote{Corresponding author. E-mail: schub.danil@yandex.ru}}
\affil{National Research University ``Higher School of Economics'' }

\maketitle
	
\begin{abstract}
In this paper we consider non-singular Morse-Smale flows on closed orientable 3-manifolds, under the assumption that among the periodic orbits of the flow there is only one saddle orbit and it is twisted. It is found that any manifold admitting such flows is either a lens space, or a connected sum of a lens space with a projective space, or Seifert manifolds with base sphere and three special layers. A complete topological classification of the described flows is obtained and the number of their equivalence classes on each admissible manifold is calculated.
\end{abstract}
	
\section{Introduction and formulation of results}
In this paper we consider, {\it NMS-flows} $f^t$, i.e., {\it not-singular} (without fixed points) Morse-Smale flows defined on closed connected orientable 3-manifolds $M^3$. The non-wandering set of such a flow consists of a finite number of periodic hyperbolic orbits. In the neighborhood of a hyperbolic periodic orbit $\mathcal O$, the flow admits a simple description (up to topological equivalence), namely, there exists its tubular neighborhood $V_{\mathcal O}$ homeomorphic to the solid torus $\mathbb V=\mathbb D^2\times \mathbb S^1$, in which the flow is topologically equivalent to the suspension over some linear diffeomorphism of the plane given by a matrix with positive determinant and real eigenvalues modulo different from unity (see e.g. \cite{Irwin}). If both eigenvalues are modulo greater than (less than) one, then the corresponding periodic orbit is {\it repelling $($attracting$)$}, otherwise it is {\it saddle}. A saddle orbit is called {\it twisted} if both eigenvalues are negative and {\it is non-twisted} otherwise.

The dynamics of NMS-flows has been studied in a number of papers: M. Wada~\cite{wada1989closed} showed that the link consisting of periodic orbits of the NMS-flow on the sphere is obtained from the Hopf link by applying a finite number of certain operations (Wada operations); J. Franks~\cite{Fr} described the flows on the three-dimensional sphere using the Lyapunov graph;
It is known from the work of Azimov \cite{Azimov} that the ambient manifold in this case is a union of circular handles. However, in the case of a small number of periodic orbits, the topology of the manifold can be substantially refined. For example, NMS-flows with exactly two periodic orbits, attracting and repelling (such a pair of periodic orbits must be contained by any NMS-flow), admit only lens spaces. 
Moreover, in \cite{PoSh} it is proved that every lens space admits exactly two equivalence classes of such flows, except for the 3-sphere $\mathbb S^3$ and the projective space $\mathbb RP^3$, on which the equivalence class is one.
In the case of a larger number of orbits, the topology of the ambient manifold is considerably richer: in the paper~\cite{Shu21} we construct NMS-flows with three periodic orbits on small Seifert manifolds.
The topology of compact orientable 3-manifolds admitting NMS-flows was studied in more detail by Morgan in~\cite{Morgan}.

In the present paper we obtain an exhaustive classification for the set $G^-_1(M^3)$ of NMS flows $f^t\colon M^3\to M^3$ with a single saddle orbit, under the assumption that it is twisted. Note that such information cannot be obtained from the general classification of Morse-Smale flows on 3-manifolds obtained in the works of Umansky~\cite{Umansky} and Prishlyak~\cite{Prishlyak}.

\begin{figure}[h!]
	\centering
	\begin{minipage}[b]{0.495\textwidth}
\center{\includegraphics[width=1\linewidth]{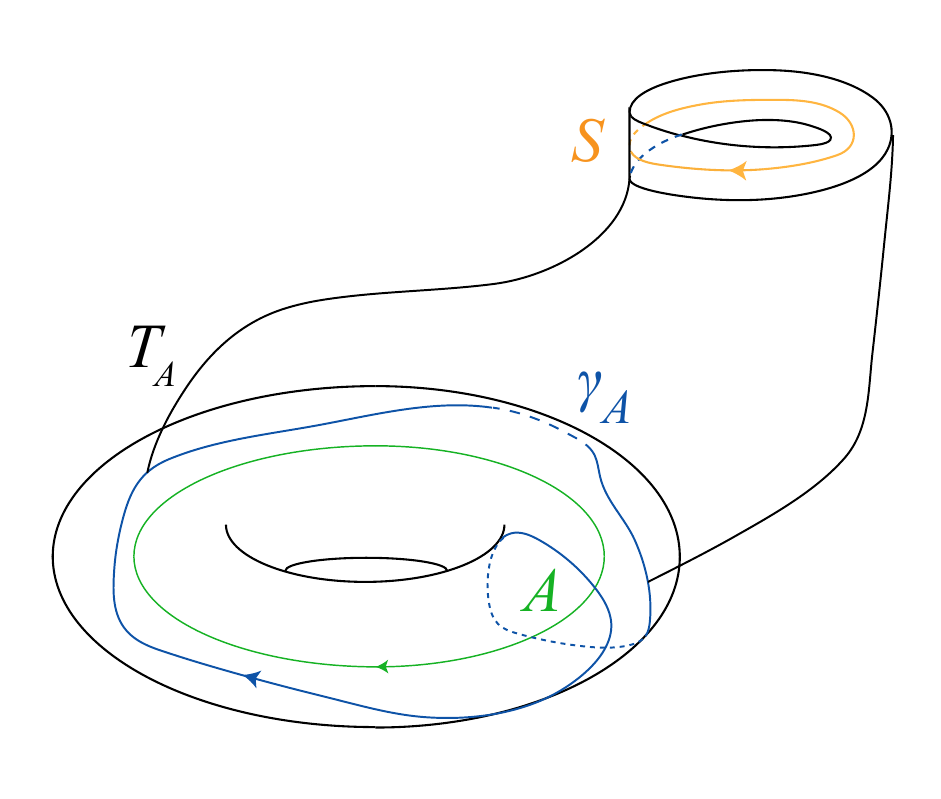}}
		\\
		{essential}
	\end{minipage}
	\hfill
	\begin{minipage}[b]{0.495\textwidth}
\center{\includegraphics[width=1\linewidth]{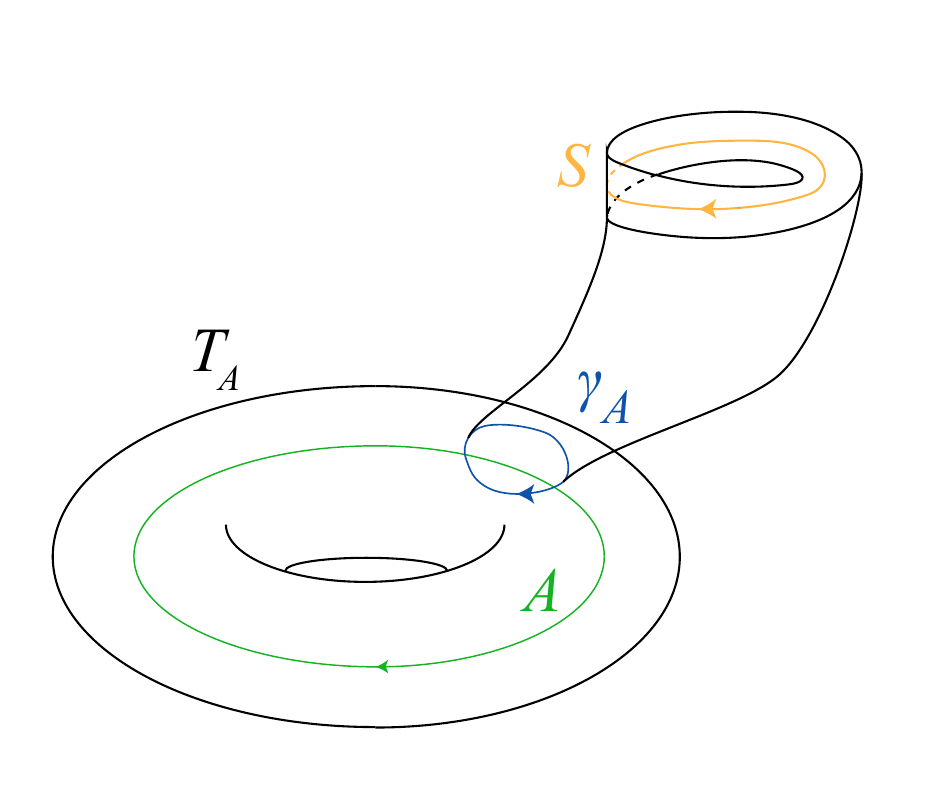}}
		\\
		{inessential}
	\end{minipage}
	\caption{Knot $\gamma_{A}$}\label{pm-}
\end{figure}

Since the ambient manifold of the flow $f^t\in G^-_1(M^3)$ is the union of stable (unstable) manifolds of all its periodic orbits \cite{Sm}, the flow must have at least one attracting and at least one repelling orbit. In the present paper the following fact is established.
\begin{lemma}\label{lem:RAS} The non-wandering set of any flow $f^t\in G^-_1(M^3)$ consists of exactly three periodic orbits $S,A,R$, saddle, attracting and repelling, respectively.
\end{lemma}
Due to the equivalence of the flow $f^t$ in the neighborhood of a periodic orbit to a suspension over a linear diffeomorphism, the unstable and stable manifolds of these orbits have the following topology:
\begin{itemize}
\item $W^u_S\cong W^s_S\cong\mathbb R\tilde\times\mathbb S^1$ (open Mobius band);
\item $W^s_A\cong W^u_R\cong\mathbb R^2\times\mathbb S^1$;
\item $W^u_A\cong W^s_R\cong\mathbb S^1$.
\end{itemize} 
A consequence of the topology of invariant manifolds of periodic orbits and the Lemma~\ref{lem:RAS} is the following representation of the $M^3$ ambient manifold.
\begin{lemma}\label{lem:union} The ambient manifold $M^3$ of any flow $f^t\in G^-_1(M^3)$ is represented as a union of three solid tori  
$$M^3 = V_A \cup V_S \cup V_R$$ with non-intersecting interior, which are tubular neighborhoods of the orbits $A,S,R$, respectively, with the following properties:
\begin{itemize}
\item torus $T_S=\partial V_S$ is the union of compact tubular neighborhoods $K_A,\,K_R$ of nodes $\gamma_{A}=W^u_S\cap T_S,\,\gamma_{R}=W^s_S\cap T_S$, respectively, such that $K_A\cap K_R=\partial K_A\cap\partial K_R$ (see Fig. \ref{pm-});
\item torus $T_A=\partial V_A$ is the union of the annulus $K_A$ and the compact surface $K=T_A\setminus \operatorname{int} K_A$,
\item torus $T_R=\partial V_R$ is the union of the annulus $K_R$ and the compact surface $K=T_R\setminus \operatorname{int} K_R$.
\end{itemize}
\end{lemma}

For $\mathcal O\in\{A,S,R\}$ we choose {\it parallel} $L_{\mathcal O}$ on the torus $T_{\mathcal O}$ (a curve homologous in $V_{\mathcal O}$ to the orbit of $\mathcal O$) and {\it meridian} $M_{\mathcal O}$ (a curve, homotopic to zero on $V_{\mathcal O}$ and essential on $T_{\mathcal O}$) such that the ordered pair of curves $L_{\mathcal O},\,M_{\mathcal O}$ defines the outer side of the solid torus $V_{\mathcal O}$.

Let $\gamma_S$ be the connected component of the set $\partial K$ oriented coherently with the saddle orbit $S$. By virtue of the equivalence of the flow $f^t|_{V_S}$ to the suspension, the meridian $M_{S}$ can be chosen such that $\gamma_S$ intersects the meridian $M_{S}$ at exactly two points. 
Then the generators $L_S,M_S$ can be chosen so that with respect to them the node $\gamma_S$ has homotopic type $$\langle \gamma_{S}\rangle=\langle l_S,m_S\rangle=\langle 2,1\rangle.$$

We orient the nodes $\gamma_{R},\, \gamma_A$ consistent with the node $\gamma_S$. Let us write down the homotopy type of the node $\gamma_R$ with respect to the $L_R,M_R$ $$\langle \gamma_{R}\rangle=\langle l_{R},m_{R}\rangle$$ and the homotopy type of the node $\gamma_A$ with respect to the $L_A,M_A$ $$\langle \gamma_{A}\rangle=\langle l_{A},m_{A}\rangle.$$
Since $T_R\setminus\gamma_R=T_A\setminus\gamma_A$, then $$(l_R, m_R)= (0,0)\iff (l_A, m_A)= (0,0).$$
If $(l_R, m_R)=(0,0)$, then let us write the homotopy type of the meridian $M_R\subset K$ with respect to the formers $L_A,\, M_A$ $$\langle M_R\rangle=\langle p_A, q_A \rangle.$$ 

If $(l_R, m_R)\neq (0, 0)$, then choose a node $\sigma_S\subset T_S$ such that
$$\langle \sigma_S\rangle = \langle 1, 1\rangle$$.
and $\sigma_S$ intersects with each component of the $\partial K$ connectivity at exactly one point (this can be done since the intersection index of the nodes $\gamma_S$ and $\sigma_S$ is $1$).

Let us choose nodes $\sigma_R\subset T_R,\,\,\sigma_A\subset T_A$ coinciding with each other on the ring $K$ and such that $\sigma_S=(\sigma_R\cup\sigma_A)\cap T_S$.
Let us write their homotopy types with respect to generators 
$$\langle \sigma_{R}\rangle=\langle b_{R},c_{R}\rangle,\,\langle \sigma_{A}\rangle=\langle b_{A},c_{A}\rangle.$$
 
\begin{definition}
	By the flow $f^t\in G^-_1(M^3)$, we define \emph{set} $$C_{f^t}=(l_1,\, b_1,\, l_2,\, b_2)$$ as follows:
	\begin{itemize}
		\item $(l_1,\, b_1,\, l_2,\, b_2)=(l_R,\, b_R,\, l_A,\, b_A)$ if $(l_R,\, m_R) \neq (0,\, 0)$;
		\item $(l_1,\, l_2,\, b_1,\, b_2)=(0,\, 2,\, p_A,\, q_A)$ if $(l_R,\, m_R) = (0,\, 0)$ and the 2-disk bounded by knot $\gamma_R$ remains on the left when moving along the knot;
		\item $(l_1,\, b_1,\, l_2,\, b_2)=(0,\, -2,\, -p_A,\, -q_A)$ if $(l_R,\, m_R) = (0,\, 0)$ and the 2-disk bounded by the knot $\gamma_R$ remains on the right when moving along the knot.
\end{itemize} 
\end{definition}

Note that the set $C_{f^t}$ of the flow $f^t\in G^-_1(M^3)$ is admissible in the sense of the following definition.

\begin{definition}
	The set of integers $C=(l_1,b_1,l_2,b_2)$ is called \emph{admissible} if
	\begin{itemize}
		\item $(l_1, b_1) = (0, \pm 2)$ or $\gcd(l_1, b_1)=1$;
		\item $\gcd(l_2, b_2)=1$.
	\end{itemize}
\end{definition}

\begin{definition}
	We call the admissible sets $C=(l_1,b_1,l_2,b_2)$, $C'=(l'_1,b'_1,l'_2,b'_2)$ \emph{consistent} $(C\sim C')$ if:
	\begin{itemize}
		\item $l_i=l'_i,\ i=1,2$, 
	\end{itemize}
	and exists $\delta\in \{-1,\, 1\}$ such that 
	\begin{itemize}
		\item $b_i \equiv \delta b'_i \pmod{l_i}$;
		\item $l_1l_2(2l_2(b_1 - \delta b'_{1})+2{l_1}(b_2 - \delta b'_{2})+l_1l_2(1 - \delta))=0$.
	\end{itemize}
\end{definition}

In the present work, the following classification result is established.

\begin{theorem}\label{th:main}
	The flows $f^t,\, f'^t\in G^-_3(M^3)$ are topologically equivalent if and only if $C_{f^t}\sim C_{f'^t}$. Moreover, for any admissible set $C$ there exists a flow $f^t\in G^-_3(M^3)$ such that $C\sim C_{f^t}$.
\end{theorem}

We also managed to construct a correspondence between invariants and ambient manifolds of flows of the considered class.

\begin{theorem}\label{th:top} 
	Flows of class $G^-_1(M^3)$ admit all lens spaces $L_{p,q}$, all connected sums of the form $L_{p,q}\#\mathbb \rp^3$ and all Seifert manifolds of the form $M(\mathbb S^2,(\alpha_1,\beta_1),(\alpha_2,\beta_2), (2, 1))$. More precisely, let the flow $f^t\in G^-_1(M^3)$ have the invariant $C_{f^t} = (l_1, b_1, l_2,b_2)$. Then
	\begin{enumerate}[label={\arabic*)}]
		\item If $(|l_1| - 1)(|l_2| - 1) = 0$, then $M^3\cong L_{p, q}$, thus:
		\begin{enumerate}[label={\roman*)}]
			\item if $l_1l_2=0$, then \\\ $M^3 \cong \rp^3$;		
			\item if $C_{f^t} = (\pm 1,\, b_1,\, l_2,\, b_2),\ l_2\neq 0$, then \\\ $M^3 \cong L_{l_2-2b_2,b_2}$;
			\item if $C_{f^t} = (l_1,\, b_1,\, \pm 1,\, b_2),\ l_1\neq 0$, then \\ $M^3 \cong L_{l_1-2b_1,b_1}$;
		\end{enumerate}
		\item If $l_1l_2 = 0$ and $(|l_1| - 1)(|l_2| - 1) \neq 0$, then $M^3\cong L_{p, q}\# \rp^3$, thus:
		\begin{enumerate}[label={\roman*)}]
			\item if $C_{f^t} = (0,\, b_1,\, l_2,\, b_2)$, then \\\ $M^3 \cong L_{l_2,\, b_2}\# \rp^3$;
			\item if $C_{f^t} = (l_1,\, b_1,\, 0,\, \pm 1),\_1 l\neq 0$, then \ $M^3 \cong L_{l_1,\, b_1}\# \rp^3$.
		\end{enumerate}
		\item If $C_{f^t} = (l_1,\, b_1,\, l_2,\, b_2),\ |l_1| > 1,\ |l_2|>1$, then $M^3\cong M(\mathbb S^2, (l_1,b_1), (l_2,b_2), (2,1)).$
\end{enumerate}
\end{theorem}

Due to the fact that the topological equivalence class and the topology of a manifold are defined using the same invariant, it becomes possible to compute the number of topological equivalence classes on each admissible manifold. 
For this purpose, for any pair $p,q$ of coprime integers, let $\bar p=|p|$ and denote by $\bar q$ -- the smallest non-negative of the numbers $q'$ satisfying the condition $q \equiv \pm q' \pmod{p}$, and by $\tilde q$ -- the smallest non-negative of the numbers $q'$ satisfying the condition $qq' \equiv \pm 1 \pmod{p}$.

\begin{theorem}\label{th:num-class} The set $G^-_1(L_{p, q}),\,|p||\neq 2$ decomposes into a countable number of equivalence classes, whereas the sets $G^-_1(\rp^3)$, $G^-_1(L_{p,q}\# \rp^3),\, G^-_1(M(\mathbb S^2,\ (\alpha_1, \beta_1),\ (\alpha_2,\beta_2),\ (2, 1)))$ consist of a finite number of classes. Namely,
\begin{enumerate}[label={\arabic*)}]
\item equivalence classes of the set $G^-_1(L_{p, q})$ depending on $p,q$ are represented by flows with the following invariants:
\begin{enumerate}
\item[a)] $|p|>2,\,q^2\not\equiv\pm 1\pmod{p}$,  $n,k\in\mathbb Z$
\begin{gather*}
\left(\pm 1, n,\bar p+2(\bar q+k\bar p), \bar q+k\bar p \right), 
\left(\pm 1,n,-\bar p+2(\bar q+k\bar p),\bar q+k\bar p \right),\\
\left(\pm 1, n,\bar p+2(-\bar q+k\bar p), -\bar q+k\bar p \right), 
\left(\pm 1,n,-\bar p+2(-\bar q+k\bar p),-\bar q+k\bar p \right),\\
\left(\pm 1,n,\bar p+2(\tilde q+k\bar p), \tilde q+k\bar p \right), 
\left(\pm 1,n,-\bar p+2(\tilde q+k\bar p), \tilde q+k\bar p \right),\\ 
\left(\pm 1,n,\bar p+2(-\tilde q+k\bar p), -\tilde q+k\bar p \right)\left(\pm 1,n,-\bar p+2(-\tilde q+k\bar p), -\tilde q+k\bar p \right),\\
\left(\bar p+2(\bar q+k\bar p), \bar q+k\bar p,\pm 1, n \right), 
\left(-\bar p+2(\bar q+k\bar p),\bar q+k\bar p,\pm 1, n, \right),\\
\left(\bar p+2(-\bar q+k\bar p), -\bar q+k\bar p,\pm 1, n \right), 
\left(-\bar p+2(-\bar q+k\bar p),-\bar q+k\bar p,\pm 1, n \right),\\
\left(\bar p+2(\tilde q+k\bar p), \tilde q+k\bar p,\pm 1, n \right), 
\left(-\bar p+2(\tilde q+k\bar p), \tilde q+k\bar p,\pm 1, n \right),\\ 
\left(\bar p+2(-\tilde q+k\bar p), -\tilde q+k\bar p,\pm 1, n \right)\left(-\bar p+2(-\tilde q+k\bar p), -\tilde q+k\bar p,\pm 1, n, \right);
\end{gather*}
\item[b)] $|p|>2,\,q^2\equiv\pm 1\pmod{p}$, $n,k\in\mathbb Z$
\begin{gather*}
\left(\pm 1, n,\bar p+2(\bar q+k\bar p), \bar q+k\bar p \right), 
\left(\pm 1,n,-\bar p+2(\bar q+k\bar p),\bar q+k\bar p \right),\\
\left(\pm 1, n,\bar p+2(-\bar q+k\bar p), -\bar q+k\bar p \right), 
\left(\pm 1,n,-\bar p+2(-\bar q+k\bar p),-\bar q+k\bar p \right),\\
\left(\bar p+2(\bar q+k\bar p), \bar q+k\bar p,\pm 1, n \right), 
\left(-\bar p+2(\bar q+k\bar p),\bar q+k\bar p,\pm 1, n, \right),\\
\left(\bar p+2(-\bar q+k\bar p), -\bar q+k\bar p,\pm 1, n \right), 
\left(-\bar p+2(-\bar q+k\bar p),-\bar q+k\bar p,\pm 1, n \right);
\end{gather*}
\item[c)] $p=0$, $n\in\mathbb Z$
\begin{gather*}
\left(\pm 1,n, 2,1 \right), 
\left(\pm1,n, -2,-1 \right),
\left(2,1,\pm 1,n \right),
\left(-2,-1, \pm 1,n \right);
	\end{gather*}	
\item[d)] $|p|=1$, $n,\,k\in\mathbb Z$
\begin{gather*}
\left(\pm 1,n, 1+2k,k \right), 
\left(1+2k,k,\pm 1,n \right);
	\end{gather*}
\item[e)] $|p|=2$ 
\begin{gather*}
\left(\pm 1, 0, 0, 1 \right), 
\left(0,1,\pm 1,0\right), 
\left(0,2,\pm 1,0 \right);
	\end{gather*}
\end{enumerate}	
\item equivalence classes of the set $G^-_1(L_{p,q}\#\rp^3)$ depending on $p,q$ are represented by flows with the following invariants:
\begin{enumerate}
\item[a)] $|p|>2,\,q^2\not\equiv\pm 1\pmod{p}$
\begin{gather*}
\left(0,2,\bar p,\pm \bar q \right), \left(0,-2,\bar p,\pm \bar q \right),\left(0,-2,\bar p,\pm \bar q \right), \left(0,-2,-\bar p,\pm \bar q \right),\\ \left(0,2,\bar p,\pm\tilde q \right),\left(0,2,-\bar p,\pm\tilde q \right),\left(0,-2,\bar p,\pm\tilde q \right),\left(0,-2,-\bar p,\pm\tilde q \right),\\
\left(0, 1,\bar p,\pm \bar q \right), \left(0, 1,-\bar p,\pm \bar q \right),\left(0, -1,\bar p,\pm \bar q \right),\left(0, -1,-\bar p,\pm \bar q \right),\\
\left(0, 1,\bar p,\pm \tilde q \right),\left(0, 1,-\bar p,\pm \tilde q \right),\left(0, -1,\bar p,\pm \tilde q \right),\left(0, -1,-\bar p,\pm \tilde q \right),\\	
\left(\bar p,\pm \bar q,0,1 \right), \left(-\bar p,\pm \bar q,0,1 \right),\left(\bar p,\pm \bar q,0,-1 \right),\left(-\bar p,\pm \bar q,0,-1 \right),\\
\left(\bar p,\pm \tilde q,0, 1 \right),\left(-\bar p,\pm \tilde q,0, 1 \right),
\left(\bar p,\pm \tilde q,0, -1 \right),
\left(-\bar p,\pm \tilde q,0, -1 \right);
\end{gather*}
\item[b)] $|p|>2,\,q^2\equiv\pm 1\pmod{p}$ 
\begin{gather*}
\left(0,2,\bar p,\pm \bar q \right), \left(0,2,-\bar p,\pm \bar q \right),\left(0,-2,\bar p,\pm \bar q \right), \left(0,-2,-\bar p,\pm \bar q \right),\\
\left(0, 1,\bar p,\pm \bar q \right), \left(0, 1,-\bar p,\pm \bar q \right),\left(0, -1,\bar p,\pm \bar q \right),\left(0, -1,-\bar p,\pm \bar q \right),\\	
\left(\bar p,\pm \bar q,0,1 \right), \left(-\bar p,\pm \bar q,0,1 \right), \left(\bar p,\pm \bar q,0,-1 \right), \left(-\bar p,\pm \bar q,0,-1 \right);
\end{gather*}
\item[c)] $p=0$
\begin{gather*}
\left(0,2, 0,\pm 1 \right), 
\left(0,1, 0,\pm 1 \right);
\end{gather*}
\item[d)] $|p|=2$
\begin{gather*}
\left(0,2, 2,\pm 1 \right), 
\left(0,2, -2,\pm 1 \right), 
\left(0,1, 2,\pm 1 \right),
\left(0,1, -2,\pm 1 \right),\\
\left(2,\pm 1,0,1 \right),
\left(-2,\pm 1,0,1 \right);
\end{gather*}
\end{enumerate}	
\item equivalence classes of the set $G^-_1(M(\mathbb S^2,\ (\alpha_1, \beta_1),\ (\alpha_2,\beta_2),\ (2, 1)))$ depending on $\alpha_1, \beta_1,\alpha_2,\beta_2$ are represented by flows with the following invariants:
\begin{enumerate}
\item[a)] $\alpha_1=\alpha_2=\alpha,\ \beta_1=\beta_2=\beta$ 
\begin{gather*} 
	\left(\alpha,\beta,\alpha,\beta \right).
\end{gather*}	
\item[b)] $|\alpha_1-\alpha_2|+|\beta_1-\beta_2|>0$ 
\begin{gather*} 
	\left(\alpha_1,\beta_1,\alpha_2,\beta_2 \right),
	\left(\alpha_2,\beta_2,\alpha_1,\beta_1 \right).
\end{gather*}
\end{enumerate}	
\end{enumerate}
\end{theorem}
Note that on the three-dimensional sphere $\mathbb S^3$, the list of flows representing equivalence classes of the set $G^-_1(\mathbb S^3)$ listed in Theorem~\ref{th:num-class} is exactly the same as that obtained in Bin Yu's paper (see ~Proposition~7.4~in~\cite{Yu}).

\textit{Acknowledgements}. The work was carried out within the framework of the fundamental research program of the National Research University Higher School of Economics.

\section{Topology of 3-manifolds}
\subsection{Lens spaces}
Futhere, we will assume that the constituents of the homotopy types of knots on the boundary $\partial\mathbb V$ of the standard fullness $\mathbb V = \mathbb D^2\times \mathbb S^1$ are the meridian $\mathbb M=(\partial \mathbb D^2) \times \mathbb S^1$ with homotopy type $\langle 0, 1\rangle$ and the parallel $\mathbb L=\{x\}\times \mathbb S^1,\,x\in \partial \mathbb D^2$ with homotopy type $\langle 1,0\rangle$.

A three-dimensional manifold $L_{p,q} = V_1 \cup_j V_2$ resulting from gluing two copies of a solid torus $V_1=\mathbb V$ is called \emph{a lens space}, $V_2= \mathbb V$ by some homeomorphism $j\colon \partial V_1\to \partial V_2$ such that $j_*(\langle 0,1\rangle)=\langle p,q\rangle$.

\begin{proposition}[\cite{Hatcher}]\label{lens-class}
Two lens spaces $L_{p,q},\,L_{p',q'}$ are homeomorphic if and only if $|p| =|p'|,\ q \equiv \pm q' \pmod{p}$ or $qq' \equiv \pm 1 \pmod{p}$. In this case, $$L_{0,1}\cong\mathbb S^2\times\mathbb S^1,\,L_{1,0}\cong\mathbb S^3,\,L_{2,1}\cong\mathbb RP^3.$$
\end{proposition}

\subsection{Dehn surgery along the knots and links}
Let 
\begin{enumerate}[label=(\alph*)]
		\item a closed 3-manifold $M$;
		\item knot $\gamma\subset M$; 
		\item tubular neighborhood $U_\gamma$ of node $\gamma$, which is a solid torus with standard generators on $\partial U_\gamma$ -- meridian $M_\gamma$ and parallel $L_\gamma$; 
		\item homeomorphism $h\colon\partial\mathbb V\to\partial U_\gamma$, which induces an isomorphism in the given generators such that  
		$h_{*}(\langle 0,1\rangle)=\langle\beta,\alpha\rangle.$
\end{enumerate}
Manifold
	$$M_{\gamma} = (M\setminus \mathrm{int }\ U_\gamma) \cup_{h}. \mathbb V$$ is called the \emph{ manifold obtained from the manifold $M$ by Dehn surgery along the knot $\gamma$ with equipment $\beta,\alpha$}.
	
Let us denote by $p_{\gamma}\colon(M \setminus \mathrm{int }\ U_\gamma) \sqcup \mathbb V\to M_{\gamma}$ the natural projection. 
Let us put $\tilde\gamma=p_{\gamma}(\{0\}\times\mathbb S^1),\,U_{\tilde\gamma}=p_{\gamma}(\mathbb V),\,\tilde h=p_{\gamma}h^{-1}\colon\partial U_\gamma\to\partial U_{\tilde\gamma}$. Then the manifold $M$ is recovered from $M_{\gamma}$ by the following inverse surgery. 
\begin{proposition}[\cite{Rolfsen}]\label{ts} Let $\gamma\subset M$ be a node with $\beta,\alpha$ and $\tilde\gamma$ be a node with $-\beta,\xi$, where $\alpha\xi\equiv 1\pmod{\beta}$. Then $$M\cong(M_{\gamma})_{\tilde\gamma}.$$ 
\end{proposition}

The Dehn surgery naturally generalizes to the case when $\gamma = \gamma_1\sqcup\dots\sqcup \gamma_r\subset M$ --- disjunctive union (link) of equipped knots. 
The resulting manifold $M_\gamma$ in this case is called \emph{the manifold obtained from the manifold $M^3$ by Dehn surgery along the equipped link $\gamma$}. 
A link $\gamma=\gamma_1\sqcup\dots\sqcup \gamma_r\subset M$ is called {\it is trivial} if the knots $\gamma_1,\dots,\gamma_r$ bound the pairwise non-overlapping 2-discs $d_1,\dots,d_r\subset M$. 
\begin{proposition}[\cite{Rolfsen}]\label{prop:dehn-connected-sum}.
	Let $\gamma = \gamma_1\sqcup\dots\sqcup \gamma_r\subset M$ be a trivial link with equipment $\beta_1,\alpha_1;\dots;\beta_r,\alpha_r$. Then
	$$M_{\gamma} \cong M\#L_{\alpha_1,\beta_1}\# \dots L_{\alpha_r,\beta_r}.$$
\end{proposition}

\subsection{Seifert fiber spaces}\label{Seif} 
A solid torus $\mathbb V$ partitioned into fibers of the form $\{x\}\times \mathbb S^1$ is called a \emph{trivially fibered solid torus}.
Consider the solid torus $\mathbb V = \mathbb D^2\times \mathbb S^1$ as a solid cylinder $\mathbb D^2\times [0, 1]$ with bases glued together by virtue of rotation by an angle $2\pi\nu/\alpha$ for coprime integers $\alpha, \nu,\ \ \alpha > 1$. The partitioning of a solid cylinder into segments of the form $\{x\}\times [0, 1]$ determines the partitioning of this solid cylinder into circles called \emph{fibers}. 
The segment $\{0\}\times [0, 1]$ generates a fiber called {\it special}, all other are ({\it non-special}) fibers of the solid torus are wrapped $\alpha$ times around the special layer and $\nu$ times around the meridian of the solid torus. The number $\alpha$ is called the \emph{multiplicity} of the singular fiber.
A solid torus with such a partition into fibers is called a \emph{nontrivially fibered solid torus} with {\it orbital invariants} $(\alpha,\nu)$.

An \emph{Seifert fiber space} --- is a compact, orientable 3-manifold $M$, partitioned into non-intersecting simple closed curves (fibers) such that each fiber has an tubular neighborhood entirely composed of fibers, layer-by-layer homeomorphic to the fibered solid torus. Such a partitioning is called \emph{Seifert fibration}. Fibers that under some such homeomorphism pass to the center of a nontrivially fibered solid torus are called \emph{special}.

The \emph{Base} of a Seifert fiber space $M$ is a compact surface $\Sigma = M/_\sim$, where $\sim$ is an equivalence relation such that $x \sim y$ if and only if $x$ and $y$ belong to the same layer.

The base of any Seyfert fiber space is a compact surface which is closed if and only if the manifold $M$ is closed; in particular, the base of any fibered solid torus is a disk (see, e.g., \cite{MatFom}). 
Thus, any Seifert fibration $M$ with a given base $\Sigma$ and orbital invariants 
$(\alpha_1,\nu_1),\dots,(\alpha_r,\nu_r),\,r\in\mathbb N$). 
is obtained from the manifold $\Sigma\times\mathbb S^1$ by Dehn surgery along the link $\gamma=\bigsqcup\limits_{i=1}^r \gamma_i$, where $\gamma_i=\{s_i\}\times\mathbb S^1,\,s_i\in\Sigma$ -- a knot with equipment $\beta_i,\alpha_i,\,\,\nu_i\beta_i\equiv 1\pmod{\alpha_i}$. Therefore, the generally accepted notation of such a Seifert fibration is as follows 
	$$M(\Sigma, (\alpha_1, \beta_1),\dots, (\alpha_r, \beta_r)).$$

Thus, {\it the orientation on the fibers} of the Seifert fibration is uniquely determined by the orientation of the circle $\mathbb S^1$ in the manifold $\Sigma\times\mathbb S^1$.

Two Seifert fibrations $M,M'$ are called \emph{isomorphic} if there exists a homeomorphism $h\colon M\to M'$ which maps the fibers of one fibration into the fibers of the other with preserving the orientation of the fibers. The homeomorphism $h$ in this case is called \emph{isomorphism of Seifert fibrations}.
It is not difficult to see that the fibrations on solid tori with orbital invariants $(\alpha, \nu)$ and $\alpha', \nu'$ are isomorphic if and only if, when $\alpha = \alpha'$ and $\nu \equiv \delta \nu' \pmod{\alpha}\ (\delta = \pm 1)$, and if $\delta = +1$, then the isomorphism preserves the orientation of the solid torus, otherwise it changes.

The following statement, which gives a criterion for isomorphism of two Seifert stratifications by their invariants, was proposed by Herbert Seifert in~\cite{seifert1933topologie}. An exposition of this statement, in notations closer to those given in this section, but only for orientation-preserving isomorphisms, can be found in the notes by Allen Hatcher~\cite{Hatcher} and the textbook by Sergey Matveev and Anatoly Fomenko~\cite{MatFom}.

\begin{proposition}\label{prop:Seifert-class}
The Seifert fibration $M(\Sigma, (\alpha_1, \beta_1),\dots, (\alpha_r, \beta_r))$ and $M'(\Sigma', (\alpha'_1, \beta'_1),\dots, (\alpha'_{r'}, \beta'_{r'}))$ are isomorphic if and only if $r =r'$ and the following conditions are satisfied for $\delta = \pm 1$ and the permutation $\sigma\colon \{1, 2, \dots, r\}\to \{1, 2, \dots, r\}$:
\begin{itemize}
	\item $\Sigma$ is homeomorphic to $\Sigma'$;
	\item $\alpha_i = \alpha'_{\sigma(i)};\,\beta_i \equiv \delta\beta'_{\sigma(i)} \pmod{\alpha_i}$ for $i\in\{1,\dots, r\}$;
	\item if the surface of $\Sigma$ is closed, then $\sum\limits_{i=1}^{r} \frac{\beta_i}{\alpha_i} = \delta \sum\limits_{i=1}^{r}\frac{\beta'_i}{\alpha'_i}$.
\end{itemize}
Thus, if $\delta = +1$, the isomorphism is orientation-preserving, and if $\delta = -1$, the isomorphism is orientation-reversing.
\end{proposition}

Note that some manifolds admit non-isomorphic Seifert fibrations.  All such manifolds are well known (see, for example~\cite{Hatcher}) and, as can be seen from the following statement, such manifolds include, for example, lens spaces.

\begin{proposition}[\cite{GaigesLange}]\label{lens-Seifert} 3-manifolds admits a Seifert fibrations with a base homeomorphic to sphere and at most two special fibers if and only if it is homeomorphic to a lens space. Thus, the list of all Seifert fibrations on lens spaces is as follows: 
	\begin{itemize}
		\item only the manifold $\mathbb S^2\times\mathbb S^1$  admits fibrations without special fibers;
		\item $M(\mathbb S^2, (\alpha, \beta))\cong L_{\beta,\alpha}$;
		\item $M(\mathbb S^2, (\alpha_1, \beta_1), (\alpha_2, \beta_2))\cong L_{p,q}$, where $p=\beta_1\alpha_2+\alpha_1\beta_2,\,q=\beta_1\nu_2+\alpha_1\xi_2$ and $\alpha_2\xi_2-\nu_2\beta_2 = 1$.
	\end{itemize}
\end{proposition}

The following statement, on the contrary, allows us to infer the non-homeomorphism of ambient manifolds from the non-isomorphism of Seifert stratifications.

\begin{proposition}[\cite{Hatcher}, Theorem 2.3]\label{prop:uniq-Seifert}
	 If two Seifert fibrations with three special fibers and a base sphere are not isomorphic, then the manifolds on which they are defined are not homeomorphic.
\end{proposition}

\section{Dynamics of flows of class $G^-_1(M^3)$}\label{dyd}

In this section we prove the lemmas given in the introduction.

Let us start with the lemma~\ref{lem:RAS}: the non-wandering set of any flow $f^t\in G^-_1(M^3)$ consists of exactly three periodic orbits $S,A,R$, saddle, attracting and repelling, respectively.

\begin{proof} The basis of the proof is the following representation of the ambient manifold $M^3$ of the NMS-flow $f^t$ with a set of periodic orbits $Per_{f^t}$ (see, e.g., \cite{Sm})  
	\begin{equation}\label{Mob}
		M^3 = \bigcup\limits_{\mathcal O \in Per_{f^t}} W^u_{\mathcal O}=\bigcup\limits_{\mathcal O \in Per_{f^t}} W^s_{\mathcal O},
	\end{equation} as well as the asymptotic behavior of invariant manifolds
	$$cl(W^u_{\mathcal O}) \setminus W^u_{\mathcal O} = \bigcup\limits_{\tilde{\mathcal O} \in Per_{f^t}\colon W^u_{\mathcal O}\cap W^s_{\mathcal O}\neq \varnothing} W^u_{\tilde{\mathcal O}},$$ $$cl(W^s_{\mathcal O}) \setminus W^s_{\mathcal O} = \bigcup\limits_{\tilde{\mathcal O} \in Per_{f^t}\colon W^s_{\mathcal O}\cap W^u_{\mathcal O}\neq \varnothing} W^s_{\tilde{\mathcal O}}.$$
	In particular, it follows from the above relations that any NMS-flow has at least one attracting orbit and at least one repelling orbit. Moreover, if the NMS-flow has a saddle periodic orbit, then the basin of any attracting orbit has a non-empty intersection with the unstable manifold of at least one saddle orbit (see ~Proposition~2.1.3~\cite{begin}) and the same situation with the basin of the repelling orbit.
	
	Let now $f^t\in G^-_1(M^3)$ and $S$ be its only saddle orbit. It follows from the relation (\ref{Mob}) that $W^u_S\setminus S$ intersects only with basins of attracting orbits. Since the set $W^u_S\setminus S$ is connected and the basins of attracting orbits are open, $W^u_S$ intersects exactly one such basin. Let us denote by $A$ the corresponding attracting orbit. Since the saddle orbit is unique, the attracting orbit is unique. Similar reasoning for $W^s_S$ leads to the existence of a single repelling orbit $R$.
\end{proof}

\subsection{Canonical neighborhoods of periodic orbits}\label{canon}
The flows admit a simple description (up to topological equivalence) in the neighborhood of a hyperbolic periodic orbit, namely, they are suspensions over linear diffeomorphisms of the plane.

Let us recall the definition of suspension. Let $\phi\colon \mathbb R^2\to \mathbb R^2$ be a diffeomorphism. 
Let us define the diffeomorphism $\Phi\colon \mathbb R^3\to \mathbb R^3$ by the formula $$\Phi(x_1,x_2,x_3) = (\phi(x_1,x_2),x_3-1). $$ 
Then the group $\{\Phi^n\}\cong\mathbb Z$ acts freely and discontinuously on $\mathbb R^3$, by virtue of which the orbit space $\Pi_\phi = \mathbb R^3/ \Phi$ is a smooth 3-manifold, and the natural projection $v_\phi\colon \mathbb R^3\to \Pi_\phi$ is a covering.
In this case, the flow $\xi^t\colon \mathbb R^3\to \mathbb R^3$ given by the formula 
$$\xi^t(x_1,x_2,x_3)=(x_1,x_2,x_3+t),$$    
induces a flow $[\phi]^t= v_\phi \xi^t v^{-1}_\phi:\Pi_\phi\to\Pi_\phi$, called {\it suspension}. 

We define the diffeomorphisms $a_{\pm1},a_2,a_0\colon \mathbb R^{2}\to \mathbb R^{2}$ by the formulas 
$$a_{\pm 1}(x_1,x_2) = (\pm 2x_1,\pm x_2/2),\,a_2(x_1,x_2)=(2x_1,2x_2),\,a_0=a_2^{-1}.$$ 
Suppose $$ V_{0}= \{ (x_1,x_2, x_3)\in \mathbb R^3 |\ 4^{x_3} x_1^2 +4^{x_3}x^2_2 \leqslant 1 \},$$  $$V_{\pm 1}= \{ (x_1,x_2, x_3)\in \mathbb R^3 |\ 4^{-x_3} x_1^2 +4^{x_3}x^2_2 \leqslant 1 \},$$ $$ V_{2}= \{ (x_1,x_2, x_3)\in \mathbb R^3 |\ 4^{-x_3} x_1^2 +4^{-x_3}x^2_2 \leqslant 1 \},$$
For $i\in\{0,-1,+1,2\}$, let $$T_i=\partial V_i,\,\mathbb V_i=v_{a_{i}}(V_i),\,\mathbb T_i=\partial{\mathbb V}_i,\,\mathbb O_i=v_{a_{i}}(Ox_3).$$

The following fact asserts canonical neighborhoods at hyperbolic periodic orbits.

\begin{proposition}[\cite{Irwin}]	For any hyperbolic periodic orbit $\mathcal O$ of a flow $f^t\colon M^3\to M^3$ defined on a closed orientable manifold $M^3$, there exists a tubular neighborhood $V_{\mathcal O}$ of the orbit $\mathcal O$ and a number $i_{\mathcal O}\in\{0, -1,+1,2\}$ such that the flow $f^t\big|_{V_{\mathcal O}}$ is topologically equivalent, via some homeomorphism $H_{\mathcal O}$, to the flow $[a_{i_{\mathcal O}}]^t|_{\mathbb V_{i_{\mathcal O}}}$.
	\label{Irw}
\end{proposition}
Let us call the neighborhood $V_{\mathcal O}=H_{\mathcal O}(\mathbb V_{i_{\mathcal O}})$ {\it canonical neighborhood} of the periodic orbit of $\mathcal O$. 

On the torus $\mathbb T_{i}$ we choose {\it longitude} $\mathbb L_{i}$ (a curve homologous in ${\mathbb V}_i$ to the orbit of $\mathbb O_i$) and {\it meridian} $\mathbb M_{i}$ (a curve, homotopic to zero on ${\mathbb V}_i$ and essential on ${\mathbb T}_i$) such that the ordered pair of curves ${\mathbb L}_i,\,{\mathbb M}_i$ defines the outer side of the solid torus ${\mathbb V}_i$.

In the proof of topological equivalence we will use the following fact, which follows from the proof of Theorem 4 and Lemma 4 in~\cite{PoSh}, and can also be found in~\cite{Umansky} (Theorem 1.1).

\begin{proposition}\label{h->H}
	The homeomorphism $h\colon \mathbb T_i\to \mathbb T_i$ for $i\in \{0, 2\}$ continues up to the homeomorphism $H\colon \mathbb V_i\to \mathbb V_i$, realizing the equivalence of the flow $[a_{i}]^t$ with itself, if and only if the induced isomorphism is of the form
	\footnote{Throughout the paper, we assume that the string $(l, m)$ is multiplied by the matrix on the left and the first element of the basis is the parallel of the torus.} $h_{*} = \begin{pmatrix}
		1 & k\\
		0 & \delta
	\end{pmatrix}$, where $\delta\in\{-1,1\},\,k\in\mathbb Z$.
\end{proposition}

The boundary of the canonical neighborhood of a saddle orbit, in contrast to an attracting or repulsing orbit, contains curves tangent to the suspension trajectories. Precisely, we denote by $\mathcal O_{x_1,x_2}$ the flow trajectory $\xi^t$ intersecting the plane $Ox_1x_2$ at a point with coordinates $(x_1,x_2,0)$. It is directly verified that the trajectory $\mathcal O_{x_1,x_2}$ intersects the surface $T_{\pm 1}$ if and only if $|x_1x_2|\leqslant\frac12$ and $(x_1,x_2)\neq(0,0)$. The trajectories touch the surface at one point if $|x_1x_2|=\frac12$, transversally intersect the surface at one point if $x_1x_2=0$, and otherwise transversally intersect the surface at two points 
$$\mathcal O_{{x_1,x_2}}\cap T_{\pm 1} = \{(x_1, x_2, x^s_3), (x_1, x_2, x_3^u)\},\ x_3^s < x_3^u.$$
\begin{figure}[h!]
	\center{\includegraphics[width = 0.7\textwidth]{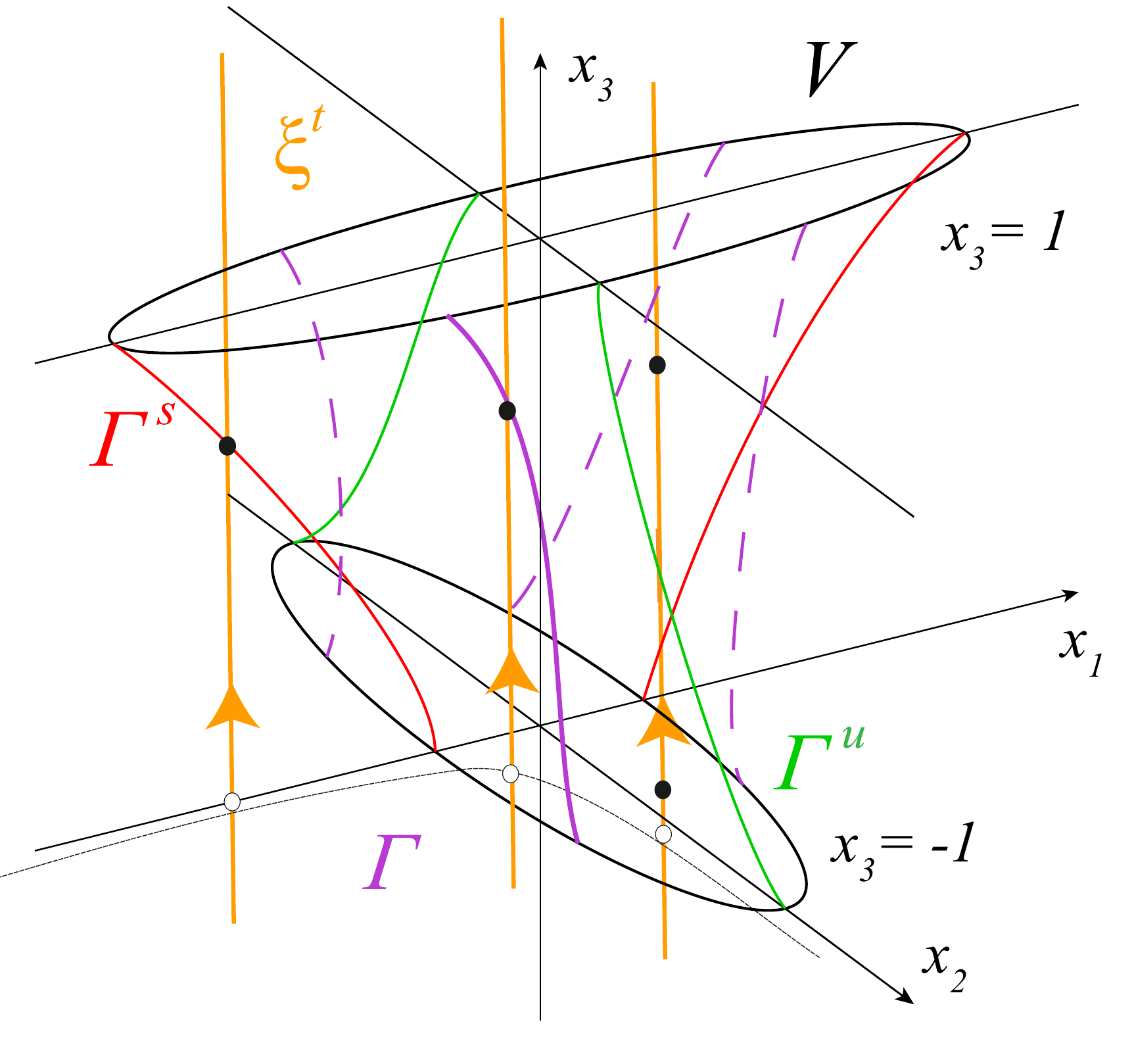}}
	\caption{Cylinder $T_{\pm 1}$}\label{V_1}
\end{figure}

Let $\Gamma=\{(x_1,x_2,x_3)\in T_{\pm 1}: \,|x_1x_2|=\frac12\}$, $\Gamma^u=Ox_1x_3\cap T_{\pm 1}$ and $\Gamma^s=Ox_2x_3\cap T_{\pm 1}$. The sets $\Gamma^u,\,\Gamma^s$ consist of two curves by construction, 
the set $\Gamma$ consist of four curves deviding $T_{\pm 1}$ into four connected components. 
The closure $T^u$ of two of these components contains $\Gamma^u$, the closure $T^s$ of two other contains $\Gamma^s$ (see.~Pic.~\ref{V_1}). We assume that $\Gamma^u$ and $\Gamma^s$ are oriented in ascending order of coordinate $x_3$. For $i\in\in\{-1,1\}$, let us put
$$\mathbb T^s_{i} = v_{a_{i}}(T^s),\ \mathbb T^u_{i} = v_{a_{i}}(T^u),\,\Gamma_{i} = v_{a_{i}}(\Gamma) ,\ \Gamma_{i}^s = v_{a_{i}}(\Gamma^s),\ \Gamma_{i}^u = v_{a_{i}}(\Gamma^u).$$

\subsection{Trajectory mappings}

In this section we prove the lemma~\ref{lem:union}:
the ambient manifold $M^3$ of any flow $f^t\in G^-_1(M^3)$ is represented as the union of three solid tori  
$$M^3 = V_A \cup V_S \cup V_R$$ with non-intersecting interior, which are tubular neighborhoods of the orbits $A,S,R$, respectively, with the following properties: 
\begin{itemize}
	\item the torus $T_S=\partial V_S$ is the union of closed tubular neighborhoods $K_A,\,K_R$ of knots $\gamma_{A}=W^u_S\cap T_S,\,\gamma_{R}=W^s_S\cap T_S$, respectively, such that $K_A\cap K_R=\partial K_A\cap\partial K_R$;
	\item the torus $T_A=\partial V_A$ is the union of the annulus $K_A$ and the compact surface $K=T_A\setminus \operatorname{int} K_A$,
	\item the torus $T_R=\partial V_R$ is the union of the ring $K_R$ and the surface $K=T_R\setminus \operatorname{int} K_R$.
\end{itemize}

\begin{proof}
Without reducing generality we will assume that the neighborhoods  $\mathcal V_A=H_A(\mathbb V_0),\ V_S=H_S(\mathbb V_{- 1}),\ \mathcal V_R=H_A(\mathbb V_2)$ of orbits $A,\,S,\,R$ pairwise disjoint.

Note, that knots $$\gamma_A = W^u_S\cap T_S=H_S(\Gamma^u_{-1}),\ \gamma_R =  W^s_S\cap T_S=H_S(\Gamma^s_{-1})$$
have tubular neighbourhoods
$$K_A =   H_S(\mathbb T^u_{-1}),\ K_R = H_S(\mathbb T^s_{- 1})$$
respectively, which are homeomorphic to annuli with a common boundary
$$\Gamma_S = H_S(\Gamma_{-1}).$$
Next, we ``blow up'' the solid tori  $\mathcal V_A$ and $\mathcal V_R$ along the trajectories so that they become  ``adjacent'' to each other and to $\mathcal V_S$. For this purpose we introduce the following notations:
\begin{itemize}
	\item let $\mathcal T_A = \partial \mathcal V_A,\quad \mathcal T_R = \partial \mathcal V_R$,
	$$K^s_R=\left(\bigcup\limits_{t>0,\,w\in cl(K_R)}f^{-t}(w)\right)\cap \mathcal T_R,\quad K^u_R = \mathcal T_R\setminus K^s_R,$$ 
	$$K^u_A=\left(\bigcup\limits_{t>0,\,w\in cl(K_A)}f^{t}(w)\right)\cap \mathcal T_A,\quad K^s_A = \mathcal T_A\setminus K^u_A;$$
	\item define continuous function $\tau_{_R}\colon \mathcal T_R\to\mathbb R^+$ such that $f^{\tau_{_R}(r)}(r)\in K_R$ for $r\in K_R$ and the set $K =\bigcup\limits_{r\in cl(K^u_R)}f^{\tau_{_{R}}(r)}(r)$ disjoint with torus $\mathcal T_A$, let $T_R=K\cup K_R$ and define homeomorphism $\psi_{_R}\colon \mathcal T_R\to T_R$ by the formula $\psi_{_R}(r)=f^{\tau_{_{R}}(r)}(r)$. Also, let $V_R$ denote the connected component of $M^3\setminus T_R$ which contains $R$;
	
	\item define continuous function $\tau_{_A}\colon \mathcal T_A\to\mathbb R^+$ such that $f^{-\tau_{_A}(a)}(a)\in K_A$ for $a\in K^u_A$ and $f^{-\tau_{_A}(a)}(a)\in K$ for $a\in K^s_A$, let $T_A=K \cup K_A$  $\psi_{A}\colon \mathcal T_A\to T_A$ and define homeomorphism $\psi_{_A}(a)=f^{-\tau_{_{A}}(a)}(a)$. Also, let $V_A$ denote the connected component of $M^3\setminus T_A$ which contains $A$;
	
	\item define continuous function $\tau_{RA}\colon K_R\setminus \gamma_R\to \mathbb R^+$ such that $f^{\tau_{RA}(w)}(w)\in K_A\setminus \gamma_A$, define homeomorphism $\psi\colon T_R\setminus \gamma_R \to T_A\setminus \gamma_A$ by the formula  
	$$\psi(w) = 
		\begin{cases}
			f^{\tau_{RA}(w)}(w), & w\in (K_R\setminus \gamma_R)\\
			w, & w\in (T_R\setminus K_R).
		\end{cases}$$
\end{itemize}

Thus, the constructed solid tori $V_A, V_S, V_R$ satisfy the conditions of the lemma.
\end{proof}

\section{Topological classification of flows $f^t\in G^-_1(M^3)$}\label{cl-}

Let us prove the first statement of Theorem \ref{th:main}: flows $f^t,\, f'^t\in G^-_1(M^3)$ are topologically equivalent if and only if their sets $C_{f^t}=(l_1,b_1,l_2,b_2),\,C_{f'^t}=(l'_1,b'_1,l'_2,b'_2)$ are consistent, that is: 
\begin{itemize}
\item $l_i=l'_i,\ i=1,2$, 
\end{itemize}
and exists $\delta\in \{-1,\, 1\}$ such that 
\begin{itemize}
\item $b_i \equiv \delta b'_i \pmod{l_i}$;
\item $l_1l_2(2l_2(b_1 - \delta b'_{1})+2{l_1}(b_2 - \delta b'_{2})+l_1l_2(1 - \delta))=0$.
		\end{itemize}
	
\begin{proof} Recall that for a periodic orbit $\mathcal O\in\{A,S,R\}$ of the flow $f^t\in G^-_1(M^3)$ we denote by $V_{\mathcal O}$ its canonical neighborhood with boundary $T_{\mathcal O}$. In this case, the ambient manifold $M^3$ of the flow $f^t$ is represented as a union of three solid tori $M^3 = V_A \cup V_S \cup V_R$ with non-intersecting interior,  torus $T_S$ is the union of compact tubular neighborhoods of $K_A,\,K_R$ knots $\gamma_{A}=W^u_S\cap T_S,\,\,\gamma_{R}=W^s_S\cap T_S$, $K=T_R\setminus{\rm int}\, K_A=T_A\setminus{\rm int}\, K_A$ and the knot $\gamma_S$ is the connected component of the boundary of the ring $K$. 

On the torus $T_{\{\mathcal O},\,\mathcal O\in\{S,R,A\}$ we have chosen the longitude $L_{\mathcal O}$ (a curve homologous in $V_{\mathcal O}$ to the orbit of $\mathcal O$) and the meridian $M_{\mathcal O}$ (a curve, homotopic to zero on $V_{\mathcal O}$ and essential on $T_{\mathcal O}$) such that the ordered pair of curves $L_{\mathcal O},\,M_{\mathcal O}$ defines the outer side of the solid torus of $V_{\mathcal O}$. 
The node $\gamma_{\mathcal O}$ is oriented coherently with the saddle orbit $S$ and has homotopy type $\langle\gamma_{\mathcal O}\rangle =\langle l_{\mathcal O},m_{\mathcal O}\rangle$ with respect to the generators $L_{\mathcal O},M_{\mathcal O}$.

If $(l_R, m_R)=(0,0)$, then we have written the homotopy type of the meridian $M_R\subset K$ with respect to the generators $L_A,\, M_A$ $$\langle M_R\rangle=\langle p_A, q_A \rangle.$$

If $(l_R, m_R) \neq (0, 0)$, then any knot $\sigma_{\mathcal O}\subset T_{\mathcal O}$ having homotopy type $\langle\sigma_{\mathcal O}\rangle =\langle b_{\mathcal O},c_{\mathcal O}\rangle$ and 
the intersection index 1 with knot $\gamma_{\mathcal O}$ has the following property  
\begin{equation}\label{mb}
l_{\mathcal O}c_{\mathcal O}
-m_{\mathcal O}b_{\mathcal O}=1.\end{equation}
Let $\Sigma^{+1}_{\mathcal O}\,(\Sigma^{-1}_{\mathcal O})$ denote set of all knots on $T_{\mathcal O}$,  having intersection index +1 (-1) with the knot $\gamma_{\mathcal O}$. Then  
\begin{equation}\label{cpm}
\tilde\sigma_{\mathcal O}\in\Sigma^{\pm1}_{\mathcal O}\iff\langle\tilde\sigma_{\mathcal O}\rangle =\langle\pm b_{\mathcal O}+n_{\mathcal O}l_{\mathcal O},\pm c_{\mathcal O}+n_{\mathcal O}m_{\mathcal O}\rangle,\,n_{\mathcal O}\in\mathbb Z.
\end{equation}
It is easily verified that the intersection index of nodes $\sigma_{\mathcal O},\,\tilde\sigma_{\mathcal O}$ is $-n_{\mathcal O}$. Then, if $\tilde\sigma_S=(\tilde\sigma_R\cup\tilde\sigma_A)\cap T_S$, then
\begin{equation}\label{nki}
n_A+n_R+n_S=0.
\end{equation}

Also recall that the generators $L_S,M_{S}$ are chosen such that with respect to them the knot $\gamma_S$ has homotopy type 
\begin{equation}\label{gaS}
\langle \gamma_{S}\rangle=\langle l_S,m_S\rangle=\langle 2,1\rangle.
\end{equation}
If $(l_R, m_R)\neq (0, 0)$ let the knot $\sigma_S\subset T_S$ be chosen such that
\begin{equation}\label{siS}
\langle \sigma_S\rangle= \langle b_S,c_S\rangle  = \langle 1, 1\rangle
\end{equation}
and $\sigma_S$ intersects with each connected component of the $\partial K$ at exactly one point (this can be done since the intersection index of the nodes $\gamma_S$ and $\sigma_S$ is $1$). Let the knots 
$\sigma_R\subset T_R,\,\sigma_A\subset T_A$ be chosen such that
\begin{equation}\label{cap}
\sigma_S=(\sigma_R\cup\sigma_A)\cap T_S.
\end{equation}
By definition
$C_{f^t}=(l_1,b_1,l_2,b_2),$ where 
\begin{itemize}
	\item $(l_1,\, b_1,\, l_2,\, b_2)=(l_R,\, b_R,\, l_A,\, b_A)$,  if $(l_R,\, m_R) \neq (0,\, 0)$;
	\item $(l_1,\, b_1,\, l_2,\, b_2)=(0,\, 2,\, p_A,\, q_A)$, if $(l_R,\,  m_R) = (0,\, 0)$ and 2-ball, bounded by the knot $\gamma_R$ remains to the left when traveling along the knot;
	\item $(l_1,\, b_1,\, l_2,\, b_2)=(0,\, -2,\, -p_A,\, -q_A)$, if $(l_R,\,  m_R) = (0,\, 0)$ and 2-ball, bounded by the knot  $\gamma_R$ remains to the right when traveling along the knot.
\end{itemize}
 
Similar equalities with primes hold for the flow $f'^t$.

Let us prove separately the necessity and sufficiency of the conditions of the theorem\ref{th:main}.

\textit{Necessity.} 
Let the flows $f^t$ and $f'^t$ with periodic orbits $A,R,S$ and $A',R',S'$ be topologically equivalent via the homeomorphism $h\colon M^3\to M^3$. For ${\mathcal O}\in\{A,S,R\}$, without reducing generality, let $V_{\mathcal O'}=h(V_{\mathcal O})$. Let $h_{\mathcal O} = h\big|_{T_{{\mathcal O}}}\colon T_{\mathcal O}\to T_{{\mathcal O}'}$.

Since $h_{\mathcal O}$ is a restriction of a homeomorphism of a solid torus, the action of the homeomorphism $h_{\mathcal O}$ in the fundamental group $\pi_1(T_{\mathcal O})$ in the generators $L_{\mathcal O},\,M_{\mathcal O}$ is given by a matrix:
\begin{equation}\label{mat}
	h_{{\mathcal O}*} = \begin{pmatrix}
		1 & k_{\mathcal O}\\
		0 & \delta_{\mathcal O}
	\end{pmatrix},\,k_{\mathcal O}\in \mathbb Z,\ \delta_{\mathcal O}\in\{-1,+1\}.
\end{equation}
Thus, since the tori $T_A,\ T_S,\ T_R$ are pairwise intersecting two-dimensional manifolds, all numbers $\delta_A,\delta_S,\delta_R$ have the same sign, let $$\delta_A=\delta_S=\delta_R=\delta_R=\delta.$$
From the properties of the conjugating homeomorphism it follows that $h_{\mathcal O}(\gamma_{\mathcal O})=\gamma_{{\mathcal O}'},\,\mathcal O\in\{S,A,R\}$, whence we obtain that  
\begin{equation}\label{(1)}
	l_{\mathcal O}=l_{\mathcal O'}
\end{equation}
and 
\begin{equation}\label{00}
(l_R,m_R)=(0,0)\iff (l_{R'},m_{R'})=(0,0).
\end{equation}
Let us prove that the consistency conditions of the sets $C_{f^t},\,C_{f'^t}$ holds separately for two cases: I) $(l_R,m_R)=(0,0)$, II) $(l_R,m_R)\neq(0,0)$.

In case I), it follows from the definition of the sets $C_{f^t},\,C_{f'^t}$ that 
$l_1=l_2=0,\,|b_2|=|b'_2|=2$. Since the homeomorphism $h_R$ maps the 2-disc bounded by the node $\gamma_R$ into the 2-disc bounded by the knot $\gamma_{R'}$ with direction of knots preserved, then $b_2=\delta b'_2$.  

It follows from the equations (\ref{mat}) that  $h_{R*}(\langle 0,1\rangle)=\langle 0,\delta\rangle$. Thus, $h_{A*}(\langle p_A,q_A\rangle)=\langle\delta p_{A'},\delta q_{A'}\rangle$. 
Since $(l_2,b_2)=(\pm p_A,\pm q_A)$, we have $(l'_2,b'_2)=(\delta (\pm p_{A'}),\delta(\pm q_{A}'))$ which implies $h_{A*}(\langle l_2,b_2\rangle)=\langle l'_2,b'_2\rangle$.  Also, it follows from the equations (\ref{mat}) that  $h_{A*}(\langle l_2,b_2\rangle)=\langle l_2,\delta b_2+k_Al_2\rangle$, whence $l_2=l_2'$ and $b_2 \equiv\delta b'_2 \pmod{l_2}$.

In case II), the equality (\ref{(1)}) is equivalent to the equality $l_i=l'_i,\,i=1,2$. Let $\tilde\sigma_{\mathcal O'}=h_{\mathcal O}(\sigma_{\mathcal O})$ and denote by $$\langle\tilde\sigma_{\mathcal O'}\rangle =\langle\tilde b_{\mathcal O'},\tilde c_{\mathcal O'}\rangle$$. 
the homotopy type of the knot $\tilde\sigma_{\mathcal O'}$ with respect to the generators $L_{\mathcal O'},M_{\mathcal O'}$. Then it follows from the formula (\ref{mat}) that  \begin{equation}\label{be}
\tilde b_{\mathcal O'}=b_{\mathcal O}
\end{equation}   
Since the determinant of the matrix $h_{\mathcal O*}$ equals $\delta$ and $h_{\mathcal O}(\gamma_{\mathcal O})=\gamma_{{\mathcal O'}}$, then $\tilde\sigma_{\mathcal O'}\in\Sigma^\delta_{\mathcal O'}$. Then from the formula (\ref{cpm}) we obtain that 
\begin{equation}\label{ura}
\tilde b_{\mathcal O'}=\delta b_{\mathcal O'}+n_{\mathcal O'}l_{\mathcal O'}, \tilde c_{\mathcal O} = \delta c_{\mathcal O'}+n_{\mathcal O'}m_{\mathcal O'}.
\end{equation}
Whence, taking into account the equalities (\ref{(1)}) and (\ref{be}), we obtain that 
\begin{equation}\label{uraa}
b_{\mathcal O}=\delta b_{\mathcal O'}+n_{\mathcal O'}l_{\mathcal O},
\end{equation}
so
\begin{equation}\label{uraaa}
	b_{\mathcal O} \equiv \delta b_{\mathcal O'} \pmod{l_{\mathcal O}}.
\end{equation} 
 
By construction $\tilde\sigma_S=(\tilde\sigma_R\cup\tilde\sigma_A)\cap T_S$, which, given the equality (\ref{nki}), entails equality  
\begin{equation}\label{nkii}
n_{A'}+n_{R'}+n_{S'}=0.
\end{equation} 
If $l_Al_R\neq 0$, then by expressing $n_{\mathcal O'}$ from the equality (\ref{uraa}) and substituting into the equality (\ref{nkii}), given that $l_S=2,\ b_S = b_{S'} = 1$, we arrive at  
$$ 2l_R(b_A - \delta b_{A'})+2{l_A}(b_R - \delta b_{R'})+l_Al_R(1 - \delta )=0,$$
which is equivalent to
$$l_Al_R(2l_R(b_A - \delta b_{A'})+2{l_A}(b_R - \delta b_{R'})+l_Al_R(1 - \delta))) = 0,$$
which is holds when $l_Al_R = 0$.

\textit{Sufficiency.} Let the sets 
$C_{f^t}=(l_1,b_1,l_2,b_2),\,C_{f'^t}=(l'_1,b'_1,l'_2,b'_2)$ of flows $f^t,\,f'^t$ are consistent via the parameter $\delta\in\{-1,1\}$. We define the homeomorphism $Q_\delta\colon\mathbb V_{-1}\to\mathbb V_{-1}$ by the formula 
$$Q_\delta = v_{a_{-1}}\bar Q_\delta v_{a_{-1}}^{-1},\text{ where } \bar Q_\delta(x_1, x_2, x_3) = (\delta x_1, x_2, x_3)\colon V_{-1}\to V_{-1}.$$

We check directly that the constructed homeomorphism $Q_\delta$ realizes the equivalence of the flow $[a_{-1}]^t$ with itself.
Let $$h_S=H_{S'}Q_\delta H^{-1}_S:V_{S}\to V_{S'}.$$

We show that the homeomorphism $h_S|_{K_A}$ can be extended to a homeomorphism $h_A\colon T_A\to T_{A'}$ inducing an isomorphism 
$h_{A*} = 
\begin{pmatrix}
	1 & k_A \\
	0 & \delta
\end{pmatrix}$ for some $k_A\in\mathbb Z$ and the homeomorphism ${h_S}|_{K_R}$ can be extended to the homeomorphism $h_R:T_R\to T_{R'}$ inducing isomorphism 
$h_{R*} = 
\begin{pmatrix}
	1 & k_R \\
	0 & \delta
\end{pmatrix}$ for some $k_R\in\mathbb Z$ such that $h_A|_{K}=h_{R}|_K$. Then, by virtue of Proposition~\ref{h->H}, the homeomorphisms $h_A,\, h_R$ can be extended to homeomorphisms $h_A\colon V_A\to V_{A'},\, h_R\colon V_R\to V_{R'}$ realizing the equivalence of the flows $f^t\big|_{V_A}$ c $ f'^t\big|_{V_{A'}}$ and $f^t\big|_{V_R}$ with $ f'^t\big|_{V_{R'}}$, respectively, and the desired 
homeomorphism $h\colon M^3\to M^3$ realizing the equivalence of $f^t,\,f'^t$ flows coincides with $h_{\mathcal O}$ on $V_{\mathcal O}$ for $\mathcal O\in\{S,A,R\}$. 

Let's consider the cases separately: I) $(l_1,b_1)=(0,\pm 2)$, II) $(l_1,b_1)\neq(0,\pm 2)$.

In case I), it follows from the consistency condition of the sets $C_{f^t},\,C_{f'^t}$ that $b_1=\delta b'_1,\,l_2=l'_2$ and $b'_2=\delta b_2+k_Al_2$ for some $k_A\in\mathbb Z$. Since the annuli $K_A,\,K_{A'}$ are contractible on tori $T_A,\,T_{A'}$, the homeomorphism $h_S|_{K_A}$ can be extended to the homeomorphism $h_A\colon T_A\to T_{A'}$ inducing an isomorphism 
$$h_{A*} = 
\begin{pmatrix}
	1 & k_A \\
	0 & \delta
\end{pmatrix}.$$
Let us define the homeomorphism $h_R\colon T_R\to T_{R'}$ by the formula
$$h_R(x) = 
\begin{cases}
	h_S(x), & x\in K_R\\
	h_A(x), & x\in K
\end{cases}.$$
Since $h_{A*}(\langle l_2,b_2\rangle)=\langle l'_2,b'_2\rangle$, then $h_{R*}(\langle 0,1\rangle)=\langle 0,\delta\rangle$
so 
$$h_{R*} = 
\begin{pmatrix}
	1 & k_R \\
	0 & \delta
\end{pmatrix}
$$ for some $k_R\in\mathbb Z$. 

In case II), it follows from the consistency condition of the sets $C_{f^t},\,C_{f'^t}$ that 
$l'_i=l_i,\,b'_{i}\equiv \delta b_{i}\pmod{l_i},\,i=1,2.$
So,
\begin{equation}\label{lb}
l_{R'}=l_R,\,b'_{R}\equiv \delta b_{R}\pmod{l_R};\,\,l_{A'}=l_A,\,b'_{A}\equiv \delta b_{A}\pmod{l_A}.
\end{equation} 
Next, we consider separately cases IIa) $l_Al_R=0$, IIb) $l_Al_R\neq 0$.

In case IIa) we assume without loss of generality that $l_R=0$ (in case $l_A=0$ the reasoning is similar). It follows from (\ref{lb}) and (\ref{mb}) that $m_{A'}=\delta m_{A}+k_Al_A$ for some $k_A\in\mathbb Z$. Then the homeomorphism $h_S|_{K_A}$ continues to a homeomorphism $h_A\colon T_A\to T_{A'}$ inducing an isomorphism 
$$h_{A*} = 
\begin{pmatrix}
	1 & k_A \\
	0 & \delta
\end{pmatrix}.$$
Let us define the homeomorphism $h_R\colon T_R\to T_{R'}$ by the formula
$$h_R(x) = 
\begin{cases}
	h_S(x), & x\in K_R\\
	h_A(x), & x\in K
\end{cases}.$$
Since $l_R=0$, then $m_R=\pm 1,\,m_{R'}=\pm\delta$ and hence $h_{R*}(\langle 0,\pm 1\rangle)=\langle 0,\pm\delta\rangle$. Then  
$$h_{R*} = 
\begin{pmatrix}
	1 & k_R \\
	0 & \delta
\end{pmatrix}
$$ for some $k_R\in\mathbb Z$.

In case IIb), the homeomorphisms $h_A\colon T_A\to T_{A'}$ and $h_R\colon T_R\to T_{R'}$ are constructed as in case IIa). Let us show that $h_{R*} =\begin{pmatrix}
	1 & k_R \\
	0 & \delta
\end{pmatrix},$ where $m_{R'}=\delta m_{R}+k_Rl_R$.

From the equality (\ref{ura}) we obtain that the knot $h_S(\sigma_S)$ has the intersection index 
$$n_{S'}=\frac{1-\delta}{2}$$ with the knot $\sigma_{S'}$, and the knot $h_A(\sigma_A)$ has the intersection index 
$$n_{A'}=\frac{b_A-\delta b_{A'}}{l_A}$$
with knot $\sigma_{A'}$.
According to the equality (\ref{nkii}), knot $h_R(\sigma_R)$ has intersection index $n_{R'}=-(n_{S'}+n_{A'})$. Then from the consistency condition of the sets, we obtain that $$n_{R'}=\frac{b_R-\delta b_{R'}}{l_R}.$$ Whence $h_{R*}(\langle b_R,c_R\rangle)=\langle b_R,\tilde c_R\rangle$. Since  $h_{R*}(\langle l_R,m_R\rangle)=\langle l_R,m_{R'}\rangle$, then 
$$h_{R*} = 
\begin{pmatrix}
	1 & k_R \\
	0 & \delta
\end{pmatrix}.
$$
\end{proof}

\section{Realization of flows $f^t\in G^-_1(M^3)$}

In this section we prove the second part of Theorem~\ref{th:main}: for any admissible invariant $C$ there exists a flow $f^t\in G^-_1(M^3)$. Recall that an invariant $C$ is called admissible if:
\begin{itemize}
	\item $(l_1, b_1) = (0,\pm 2)$ or $\gcd(l_1, b_1)=1$;
	\item $\gcd(l_2, b_2)=1$.
\end{itemize}

\begin{proof}
	Let $C = (l_1,\, b_1,\ l_1,\ l_1,\, b_1)$. 
	We construct the three-dimensional manifold $M^3$ and the flow $f^t\in G^-_1(M^3)$ such that $C_{f^t} = C$ separately for the cases: I) $(l_1, b_1) = (0,\pm 2)$, II) $(l_1, b_1) \neq(0,\pm 2)$.
	
	In case I), for $(l_1,\, b_1) = (0, \pm 2)$, let $(p,q) = (\pm l_2,\pm b_2)$. Let us define a homeomorphism $\psi\colon\mathbb T_2\to\mathbb T_0$ inducing an isomorphism defined by the integer matrix $h_{*} = 
		\begin{pmatrix}
			r & s \\
			p & q
		\end{pmatrix}$ 
		with determinant equals to -1. On the torus $\mathbb T_2$, we choose a $\gamma_R$ essential knot with a tubular neighborhood $K_R\subset (\mathbb T_2\setminus\mathbb M_2)$ and orient it so that the 2-disk bounded by it remains on the left in the case $b_1=+2$ and -- on the right in the case $b_1=-2$. Let $\gamma_A=\psi(\gamma_R)$ and $K_A=\psi(K_R)$.

		Let $\psi_{R}\colon K_R\to \mathbb T^s_{-1},\\psi_{A}\colon K_A\to \mathbb T^u_{-1}$ be such homeomorphisms, that $\psi^{-1}_A\psi_R|_{\partial K_R}=\psi|_{\partial K_R}$, $\psi_{R}(\gamma_R)=\Gamma^s$, $\psi_{A}(\gamma_A)=\Gamma^u$. Let $\sim$ be the minimal equivalence relation on $\mathbb V_0 \sqcup \mathbb V_{-1} \sqcup \mathbb V_2$ for which $x\sim \psi(x),x\in (\mathbb T_2\setminus {\rm int}\,K_A),\ x\sim \psi_{A}(x),x\in K_A,\ x\sim \psi_{R}(x),x\in K_R$. Then 
		$$M^3 = (\mathbb V_0 \sqcup \mathbb V_{-1} \sqcup \mathbb V_2)/\sim.$$
		We denote by $\pi\colon \mathbb V_0 \sqcup \mathbb V_{-1} \sqcup \mathbb V_2 \to M^3$ the natural projection.
		Let the flow $f^t\colon M^3\to M^3$ be given by the formula
		$$f^t(x) = \begin{cases}
			\pi([a_0]^t(\pi^{-1}(x))), & x\in \pi(\mathbb V_0)\\
			\pi([a_{-1}]^t(\pi^{-1}(x))), & x\in \pi(\mathbb V_{-1})\\
			\pi([a_2]^t(\pi^{-1}(x)), & x\in \pi(\mathbb V_2))
		\end{cases}$$
		By construction, $C_{f^t} = C$.
	
		In case II), we represent the sphere $\mathbb S^2$ as a union of three two-dimensional disks $D_{A},\ D_{S},\ D_{R}$ with centers $O_A,\ O_S,\ O_R$, glued along the boundary as depicted in Picture~\ref{base2} (glued segments are marked with the same color).
		\begin{figure}[!h]
		\center{\includegraphics[width=0.7\textwidth]{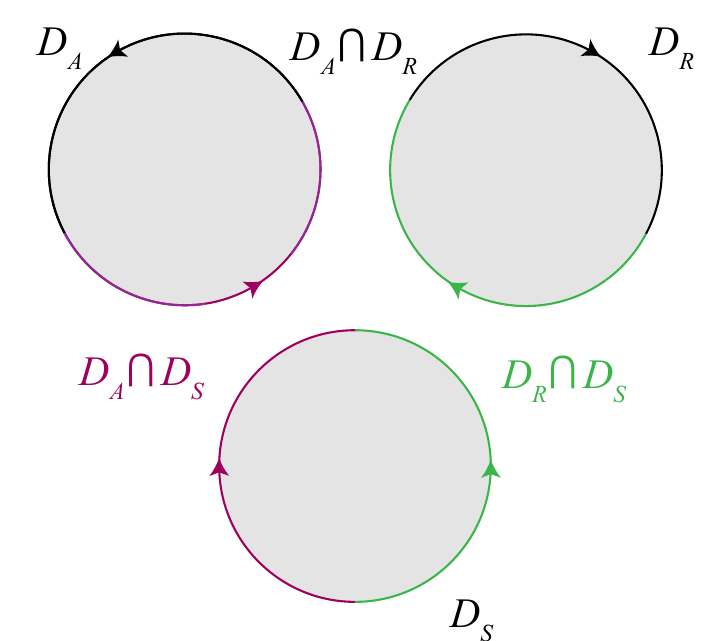}}
		\caption{Disks $D_A$, $D_S$, $D_R$}\label{base2}
\end{figure}
		Then the manifold $\mathbb S^2\times\mathbb S^1$ is represented as a union of three solid tori $V_A=D_A\times\mathbb S^1$, $V_S=D_S\times\mathbb S^1$, $V_R=D_R\times\mathbb S^1$, which are tubular neighborhoods of the knots $\ell_A=O_A\times\mathbb S^1$, $\ell_S=O_S\times\mathbb S^1$, $\ell_R=O_R\times\mathbb S^1$, glued along the boundaries $T_A=\partial V_A$, $T_S=\partial V_S$, $T_R=\partial V_R$, along the annuli $K_A=T_A\cap T_S$, $K_R=T_R\cap T_S$, $K=T_A\cap T_R$, $K=T_A\cap T_R$, respectively. 

		Let $C=(l_1, b_1, l_2, b_2)$.
		Then the sought manifold $M^3$ is obtained by the Dehn surgery along the link $\ell_A\sqcup\ell_S\sqcup\ell_R$ with equipment $(-b_2,l_2),(-1,2),(-b_1,l_1)$. Moreover, the homeomorphisms of surgery $h_A\colon\mathbb V_0\to V_A,\,h_S\colon\mathbb V_{-1}\to V_S,\,h_R:\mathbb V_{-1}\to V_S,\,h_R\colon \mathbb V_2\to V_R$ are chosen such that $h_S(\mathbb T^u_{-1})=K_A,\ h_S(\mathbb T^s_{-1})=K_R$. 
		Denote by $\pi\colon  \mathbb V_0 \sqcup \mathbb V_{-1} \sqcup \mathbb V_2 \to M^3$ the natural projection.
		Let $f^t\colon M^3\to M^3$ be defined by the formula 
		$$f^t(x) = \begin{cases}
			\pi([a_0]^t(\pi^{-1}(x))), & x\in \pi(\mathbb V_0)\\
			\pi([a_{-1}]^t(\pi^{-1}(x))), & x\in \pi(\mathbb V_{-1})\\
			\pi([a_2]^t(\pi^{-1}(x))), & x\in \pi(\mathbb V_2)
		\end{cases}$$
		By construction $C_{f^t} = C$.
\end{proof}

\section{Topology of ambient manifolds of flows  $f^t\in G^-_1(M^3)$}
In this section we prove the theorem~\ref{th:top}: flows of class $G^-_1(M^3)$ admit all lens spaces $L_{p,q}$, all connected sums of the form $L_{p,q}\#\mathbb \rp^3$, and all Seifert fiber spaces of the form $M(\mathbb S^2,(\alpha_1,\beta_1),(\alpha_2,\beta_2), (2, 1))$. Exactly, let the flow $f^t\in G^-_1(M^3)$ have the invariant $C_{f^t} = (l_1, b_1, l_2,b_2)$. Then
\begin{enumerate}[label={\arabic*)}]
	\item If $(|l_1| - 1)(|l_2| - 1) = 0$, then $M^3\cong L_{p, q}$, herewith:
	\begin{enumerate}[label={\roman*)}]
		\item if $l_1l_2=0$, then \\  $M^3 \cong \rp^3$;		
		\item if $C_{f^t} = (\pm 1,\, b_1,\, l_2,\, b_2),\ l_2\neq 0$, then \\ $M^3 \cong L_{l_2-2b_2,b_2}$;
		\item if $C_{f^t} = (l_1,\, b_1,\, \pm 1,\, b_2),\ l_1\neq 0$, then \\  $M^3 \cong L_{l_1-2b_1,b_1}$;
		
	\end{enumerate}
	\item If $l_1l_2 = 0$ and $(|l_1| - 1)(|l_2| - 1) \neq 0$, then  $M^3\cong L_{p, q}\# \rp^3$, herewith:
	\begin{enumerate}[label={\roman*)}]
		\item if $C_{f^t} = (0,\, b_1,\, l_2,\, b_2)$, then \\ $M^3 \cong L_{l_2,\, b_2}\# \rp^3$;
		\item if $C_{f^t} = (l_1,\, b_1,\, 0,\, \pm 1),\ l_1\neq 0$, then \\ $M^3 \cong L_{l_1,\, b_1}\# \rp^3$.
	\end{enumerate}
	\item If $C_{f^t} = (l_1,\, b_1,\, l_2,\, b_2),\ |l_1| > 1,\ |l_2|>1$, then \\ $M^3\cong M(\mathbb S^2,(l_1,b_1),(l_2,b_2),(2,1))$.
\end{enumerate}

\begin{proof} Let us prove the theorem separately for the cases: I) $l_1l_2 = 0$, II) $l_1l_2 \neq 0$. 

In case I), we denote by $M^3_S$ the manifold obtained by Dehn surgery along knot $S$ with equipment $(1,1)$ in the generators $L_S,M_S$. Let $v_S\colon (M^3 \setminus \operatorname{int }\ V_S) \sqcup \mathbb V\to M^3_S$ --- the natural projection. For simplicity, we keep the labels of all objects on $v_S(M^3\setminus \operatorname{int } V_S)$ the same as they were on $M^3\setminus \operatorname{int } V_S$ and put $\tilde S=v_S(\{0\}\times\mathbb S^1),\,V_{\tilde S} = v_S(\mathbb V)$. Then $\tilde{V}_R= V_R\cup V_{\tilde S}$ is a solid torus with boundary $\tilde T_R$ and there exists an isotopy $\zeta_t\colon V_R\to\tilde{V}_R,\,t\in[0,1]$ such that $\zeta_0=id|_{V_A},\,\zeta_t|_{K}=id|_{K},\,t\in[0,1],\,\zeta_1(V_R)=\tilde{V}_R,\,\zeta_1(\sigma_S\cap K_R)=\sigma_S\cap K_A$. For any curve $c\subset T_R$, let us put $\tilde c=\zeta_1(c)\subset\tilde T_R$. Then the isomorphism $\zeta_{1*}$ is identical in the generators $L_R,M_R;\,\tilde L_R,\tilde M_R$ and  
\begin{equation}\label{pq}
M^3_S=\tilde{V}_R\cup_\psi V_A,
\end{equation} 
where $\psi\colon\partial\tilde{V}_R\to\partial{V}_A$ is homeomorphism inducing isomorphism in generators $\tilde L_R,\tilde M_R;\,L_A,M_A$  and $\psi_{*} = 
\begin{pmatrix}
	r & s \\
	p & q
\end{pmatrix}.$ 
Hence, $M^3_S \cong L_{p, q}$. From the Statement~\ref{ts}, we obtain that 
$M^3\cong(L_{p,q})_{\tilde S},$ where $\tilde S$ is a knot with equipment $(-1,2)$. Since the knot $\tilde S$ bounds a 2-ball on at least one of the tori $\tilde{V}_A, V_R$, then, by virtue of the Statement~\ref{prop:dehn-connected-sum} 
$(L_{p,q})_{\tilde S}\cong L_{p,q}\#L_{2,-1}.$ Whence, by virtue of the Statement~\ref{lens-class}, \begin{equation}\label{dehn1}
M^3\cong L_{p,q}\#\mathbb RP^3.
\end{equation}
Let us show how the proof of 2) and 1i) follows from the deductions made.

1i)+2i) $(l_1,b_1)=(0,\pm 2)$. It follows from the definition of the set $C_{f^t}$ that $\psi_*(\langle0,1\rangle)=\langle\pm l_2,\pm b_2\rangle$. By virtue of Statement~\ref{prop:dehn-connected-sum}, $${\rm 2i)}\,M^3\cong L_{l_2,b_2}\#\rp^3,\,|l_2|\neq 1;$$ $${\rm 1i)}\,M^3\cong \rp^3,\,|l_2|= 1.$$ 
 
1i+2i) $(l_1,b_1)\neq(0,\pm 2)$.  It follows from the definition of the set $C_{f^t}$ that $\psi_*(\langle0,1\rangle)=\langle\pm l_2,\pm m_A\rangle,\,\psi_*(\langle 1,\pm c_R\rangle)=\langle\pm b_2,\pm c_A\rangle$. By direct calculation we obtain that $\psi^{-1}_*(\langle0,1\rangle)=\langle\pm l,b\rangle,$ where $|l|=|l_2|,\,|b|\equiv|b_2|\pmod{l_2}$. By virtue of Statement~\ref{prop:dehn-connected-sum}, $${\rm 2i)}\,M^3\cong L_{l_2,b_2}\#\rp^3,\,|l_2|\neq 1;$$ $${\rm 1i)}\,M^3\cong \rp^3,\,|l_2|= 1.$$

1i+2ii). By reasoning analogous to the above, we obtain that $${\rm 2ii)}\,M^3\cong L_{l_1,b_1}\#\rp^3,\,|l_1|\neq 1;$$ $${\rm 1i)}\,M^3\cong \rp^3,\,|l_1|= 1.$$

In case II), consider first the subcase $|l_1|=1$. Then $\tilde V_S=V_S\cup V_R$ is a filled torus with boundary $\tilde T_S$ and there exists an isotopy $\zeta_t\colon V_S\to\tilde{V}_S,\,t\in[0,1]$ such that $\zeta_0=id|_{V_S},\,\zeta_t|_{K_A}=id|_{K_A},\,t\in[0,1],\,\zeta_1(V_S)=\tilde{V}_S,\,\zeta_1(\sigma_S\cap K_R)=\sigma_R\cap K$. For any curve $c\subset T_S$, let $\tilde c=\zeta_1(c)\subset\tilde T_S$. Then the isomorphism $\zeta_{1*}$ is identical in the generators $L_S,M_S;\,\tilde L_S,\tilde M_S$ and  
\begin{equation}\label{pqr}
M^3_S=\tilde{V}_S\cup_\psi V_A,
\end{equation} 
where $\psi\colon\partial\tilde{V}_R\to\partial{V}_A$ is a homeomorphism, inducing in the generators $\tilde L_S,\tilde M_S;\,L_A,M_A$ isomorphism $\psi_{*} = 
\begin{pmatrix}
	r & s \\
	p & q
\end{pmatrix}.$ 
Hence, $M^3_S \cong L_{p, q}$. From the definition of the set $C_{f^t}$, it follows that $\psi_*(\langle2,1\rangle)=\langle\pm l_2, m_A\rangle,\,\,\psi_*(\langle 1,1\rangle)=\langle\pm b_2,\pm c_A\rangle$. By direct calculation, we obtain that $\psi^{-1}_*(\langle0,1\rangle)=\langle\pm l,b\rangle,$ where $|l|=||l_2-2b_2|,\,|b|\equiv|b_2|\pmod{l_2}$.  By virtue of the Statement~\ref{prop:dehn-connected-sum}, $${\rm 1ii)}\,M^3\cong L_{l_2-2b_2,b_2}.$$

In the case $|l_2|=1$ by similar reasoning, we obtain that $${\rm 1iii)}\,M^3\cong L_{l_1-2b_1,b_1}.$$

In the case $|l_1| > 1,\ |l_2| > 1$, it follows from the procedure for realizing a flow over an admissible set (see the proof of the second part of Theorem~\ref{th:main} in Case II)) that $M^3$ is a Seifert fiber space with base sphere with three special fibers 
$$M^3 \cong M(\mathbb S^2, (l_1, b_1), (l_2, b_2), (2, 1)).$$
\end{proof}

\section{Counting the number of topological equivalence classes}
In this section we give a proof of Theorem~\ref{th:num-class}. To do this, recall that for any pair $p,q$ of integer prime numbers, we put $\bar p=|p|$ and denote by $\bar q$ the smallest non-negative of the numbers $q'$ satisfying the condition $q \equiv \pm q' \pmod{p}$, and by $\tilde q$ the smallest non-negative of the numbers $q'$ satisfying the condition $qq' \equiv \pm 1 \pmod{p}$.  

{\bf Theorem~\ref{th:num-class}.} {\it The set $G^-_1(L_{p, q}),\,|p|\neq 2$ decomposes into a countable number of equivalence classes, whereas the sets $G^-_1(\rp^3)$, $G^-_1(L_{p,q}\# \rp^3),\, G^-_1(M(\mathbb S^2,\ (\alpha_1, \beta_1),\ (\alpha_2,\beta_2),\ (2, 1)))$ consist of a finite number of classes.}
\begin{proof} By virtue of Theorem~\ref{th:top}, flows of class $G^-_1(M^3)$ admit three types of manifolds 1) $L_{p, q}$; 2) $ L_{p, q} \# \rp^3$; 3) $M(\mathbb S^2,(l_1,b_1),(l_2,b_2),(2,1))$. 
	Let us prove the proof separately for each of these cases. 

	1) According to Statement~\ref{lens-class}, two lens spaces $L_{p,q},\,L_{p',q'}$ 
	are homeomorphic if and only if $\bar p =\bar p'$ and either $\bar q =\bar q'$ or $\bar q=\tilde q'$. Whence it follows that 
	$L_{p,q}\cong L_{\bar p,\bar q}$ and, $L_{p',q'}\cong L_{\bar p,\bar q}$ if and only if at least one of the following conditions for $k\in\mathbb Z$ is satisfied:
\begin{equation}\label{p+q+}
p'=\bar p,\,q'=\bar q+k\bar p;
\end{equation}
\begin{equation}\label{p-q+}
p'=-\bar p,\,q'=\bar q+k\bar p;
\end{equation}
\begin{equation}\label{p+q-}
p'=\bar p,\,q'=-\bar q+k\bar {p};
\end{equation}
\begin{equation}\label{p-q-}
p'=-\bar p,\,q'=-\bar q+k\bar {p};
\end{equation}
\begin{equation}\label{p+qq+}
p'=\bar p,\,q'=\tilde q+k\bar p;
\end{equation}
\begin{equation}\label{p-qq+}
p'=-\bar p,\,q'=\tilde q+k\bar p;
\end{equation}
\begin{equation}\label{p+qq-}
p'=\bar p,\,q'=-\tilde q+k\bar {p};
\end{equation}
\begin{equation}\label{p-qq-}
p'=-\bar p,\,q'=-\tilde q+k\bar {p};
\end{equation}
By virtue of the \ref{th:top} theorem, the lens $L_{\bar p,\bar q},\,\bar p\neq 2$ is a ambient manifold for flows with invariants  
\begin{equation}\label{inv1}
(\pm 1,n,p'+2q',q');
\end{equation} 
\begin{equation}\label{inv2}
(p'+2q',q',\pm 1,n),
\end{equation}
where $n\in\mathbb Z$. Substituting the condition (\ref{p+q+}) into (\ref{inv1}), we obtain sets of the form $$(\pm 1,b_1,\bar p+2(\bar q+k\bar p),\bar q+k\bar p).$$
From the definition of consistency, it follows that two sets of 
$$(\pm 1,n_1,\bar p+2(\bar q+k_1\bar p),\bar q+k_1\bar p),(\pm 1,n_2,\bar p+2(\bar q+k_2\bar p),\bar q+k_2\bar p),\bar q+k_2\bar p)$$ 
are consistent if and only if $k_1=k_2,\,n_1=n_2$. Thus, each representation of the lens $L_{p,q}$ in the form (\ref{p+q+}) gives rise to the family $\left(\pm 1, n,\bar p+2(\bar q+k+k\bar p), \bar q+k\bar p \right),\,n,k\in\mathbb Z$ of pairwise non-consistent sets corresponding, by virtue of Theorem~\ref{th:top}, to pairwise non-equivalent flows. If $|p|>2$, then similar families are obtained from each of the representations (\ref{p-q+}), (\ref{p+q-}), (\ref{p-q-}), (\ref{p-q-}). It is directly verified that the sets of all four families are not pairwise equivalent. Finally, if $\bar q\neq\tilde q$ (equivalent to $q^2\not\equiv\pm 1\pmod{p}$), we obtain four more families of pairwise non-equivalent sets corresponding to the representations (\ref{p+qq+}), (\ref{p-qq+}), (\ref{p-qqq+}), (\ref{p+qq-}), (\ref{p+qq-}), (\ref{p-qqq-}). Adding to the list of sets, sets of type (\ref{inv2}), we obtain a list of eight more pairwise non-equivalent sets, directly from which follows the result of the theorem in cases 1a), 1b). 

In cases 1c); 1d), by directly substituting the pairs $\bar p=0,\bar q=1$; $\bar p=1,\bar q=0$ into the sets 1b), respectively, we obtain the announced lists of pairwise non-equivalent pairs. 

In the case $|p|=2$, the lens $L_{p,q}$ is a ambient manifold for the flows with invariants  
\begin{equation}\label{inv11}
(0,c,\pm 1,n);
\end{equation} 
\begin{equation}\label{inv22}
(\pm 1,n,0,d),
\end{equation}
where $n\in\mathbb Z,\,c\in\{-2,-1,1,2\},\,d\in\{-1,1\}$. From the definition of consistency, it follows that the two sets 
$$(0,c_1,\pm 1,n_1),(0,c_2,\pm 1,n_2)$$ 
are consistent if and only if $c_1=\pm c_2,\,n_1\equiv \pm n_2\pmod 1$. A similar statement is true for sets of the form (\ref{inv22}), resulting in the announced list 1e). 

2) By virtue of Theorem~\ref{th:top} the manifold $L_{p,q}\#\rp^3,\, |p||\neq 1$ is an ambient manifold for flows with invariants  
\begin{equation}\label{inv111}
(0,c,p',q');
\end{equation} 
\begin{equation}\label{inv222}
(p',q',0,d),
\end{equation}
where $n\in\mathbb Z,\,c\in\{-2,-1,1,2\},\,d\in\{-1,1\}$. It follows from the definition of consistency that for $|p|>2$ the two sets of 
$$(0,c_1,\bar p,\bar q+k_1\bar p),(0,c_2,\bar p,\bar q+k_2\bar p)$$
are consistent if and only if $c_1=c_2,\,k_1\equiv k_2\pmod 1$. Thus, each representation of the lens $L_{p,q}$ in the form (\ref{p+q+}) gives rise to a family $(0,c,\bar p,\bar q)$ of pairwise non-consistent sets corresponding, by virtue of Theorem~\ref{th:top}, to pairwise non-equivalent flows. Similar families are obtained from each of the representations (\ref{p-q+}), (\ref{p+q-}), (\ref{p+q-}), (\ref{p-q-}), (\ref{p+qq+}), (\ref{p-qqq+}), (\ref{p+qq-}), (\ref{p-qq-}), (\ref{p-qq-}) if $\bar q\neq\tilde q$. Adding to the list of sets, the sets of type (\ref{inv222}), we obtain the list of sets announced in (2a), 2b) of this theorem.

In cases 2c); 2d), by directly substituting into sets 2b) the pairs $\bar p=0,\bar q=1$; $\bar p=2,\bar q=1$, respectively, we obtain the announced lists of pairwise non-equivalent pairs.

3) By virtue of the \ref{th:top} theorem, the manifold $M(\mathbb S^2,\ (\alpha_1, \beta_1),\ (\alpha_2,\beta_2),\ (2, 1))$ is an ambient manifold for flows with invariants  
\begin{equation}\label{inv}
(\alpha_1, \beta_1,\alpha_2,\beta_2).
\end{equation} 
By virtue of Statement~\ref{prop:Seifert-class} 
$M(\mathbb S^2,\ (\alpha_1, \beta_1),\ (\alpha_2,\beta_2),\ (2, 1))\cong M(\mathbb S^2,\ (\alpha'_1, \beta'_1),\ (\alpha'_2,\beta'_2),\ (2, 1))$
if and only if at least one of the following conditions is met:
\begin{equation}\label{to+}
\alpha'_1=\alpha_1,\,\beta'_1=\beta_1+k_1\alpha_1,\,\alpha'_2=\alpha_2,\,\beta'_2=\beta_2+k_2\alpha_2;
\end{equation}
\begin{equation}\label{to-}
\alpha'_1=\alpha_2,\,\beta'_1=\beta_2+k_1\alpha_2,\,\alpha'_2=\alpha_1,\,\beta'_2=\beta_1+k_2\alpha_1,
\end{equation}
where $k_1,k_2\in\mathbb Z$. By virtue of Theorem~ref{th:top}. 
$M(\mathbb S^2,\ (\alpha_1, \beta_1),\ (\alpha_2,\beta_2),\ (2, 1))$ is an ambient manifold for flows with invariants 
\begin{equation}\label{ab}
(\alpha'_1, \beta'_1,\alpha'_2,\beta'_2).
\end{equation}
Substituting (\ref{to+}) into (\ref{ab}), we obtain sets of the form 
$$(\alpha_1,\beta_1+k_1\alpha_1,\alpha_2,\beta_2+k_2\alpha_2).$$
It follows from the definition of consistency that all such sets are equivalent to the set $(\alpha_1, \beta_1,\alpha_2,\beta_2)$. A similar situation is obtained with the relation (\ref{to-}). From where we obtain the sets announced in 3a), 3b).
\end{proof}

\bibliographystyle{plain} 
\bibliography{refs}

\begin{thebibliography}{10}

\bibitem{Azimov}
Daniel Asimov.
\newblock Round handles and non-singular morse-smale flows.
\newblock {\em Annals of Mathematics}, 102(1):41--54, 1975.

\bibitem{Fr}
John Franks.
\newblock Nonsingular smale flows on $s^3$.
\newblock {\em Topology}, 24(3):265--282, 1985.

\bibitem{GaigesLange}
Hansj{\"o}rg Geiges and Christian Lange.
\newblock Seifert fibrations of lens spaces.
\newblock In {\em Abhandlungen aus dem Mathematischen Seminar der
  Universit{\"a}t Hamburg}, volume~88, pages 1--22. Springer, 2018.

\bibitem{begin}
Viacheslav~Z Grines, Timur~V Medvedev, and Olga~V Pochinka.
\newblock {\em Dynamical systems on 2-and 3-manifolds}, volume~46.
\newblock Springer, 2016.

\bibitem{Hatcher}
Allen Hatcher.
\newblock Notes on basic 3-manifold topology, 2007.

\bibitem{Irwin}
M.C. Irwin.
\newblock A classification of elementary cycles.
\newblock {\em Topology}, 9(1):35--47, 1970.

\bibitem{MatFom}
Sergey~Vladimirovich Matveev and Anatoly~Timofeevich Fomenko.
\newblock {\em Algorithmic and computer methods in three-dimensional topology}.
\newblock Federal State Budgetary Educational Institution of Higher Education
  "Lomonosov Moscow State University" Publishing House (Printing House), 1991.

\bibitem{Morgan}
John~W Morgan.
\newblock Non-singular morse--smale flows on 3-dimensional manifolds.
\newblock {\em Topology}, 18(1):41--53, 1979.

\bibitem{PoSh}
O.V. Pochinka and D.D. Shubin.
\newblock Non-singular morse--smale flows on n-manifolds with
  attractor--repeller dynamics.
\newblock {\em Nonlinearity}, 35(3):1485, 2022.

\bibitem{Prishlyak}
Alexander~Olegovich Prishlyak.
\newblock Complete topological invariant of morse-smele flows and handle
  decompositions of three-dimensional manifolds (in russian).
\newblock {\em Fundamental and applied mathematics}, 11(4):185--196, 2005.

\bibitem{Rolfsen}
Dale Rolfsen.
\newblock {\em Knots and links}, volume 346.
\newblock American Mathematical Soc., 2003.

\bibitem{seifert1933topologie}
Herbert Seifert.
\newblock Topologie dreidimensionaler gefaserter r{\"a}ume.
\newblock 1933.

\bibitem{Shu21}
Danila~Denisovich Shubin.
\newblock Topology of ambient manifolds of nonsingular flows with three
  nontwisted orbits.
\newblock {\em Izvestiya VUZ. Applied Nonlinear Dynamics}, 29(6):863--868,
  2021.

\bibitem{Sm}
Stephen Smale.
\newblock Differentiable dynamical systems.
\newblock {\em Bulletin of the American mathematical Society}, 73(6):747--817,
  1967.

\bibitem{Umansky}
Ya.~L. Umanskii.
\newblock Necessary and sufficient conditions for topological equivalence of
  three-dimensional morse--smale dynamical systems with a finite number of
  singular trajectories.
\newblock {\em Matematicheskii Sbornik}, 181(2):212--239, 1990.

\bibitem{wada1989closed}
Masaaki Wada.
\newblock Closed orbits of non-singular morse-smale flows on s3.
\newblock {\em Journal of the Mathematical Society of Japan}, 41(3):405--413,
  1989.

\bibitem{Yu}
Bin Yu.
\newblock Behavior $0$ nonsingular morse smale flows on $s^3$.
\newblock {\em Discrete \& Continuous Dynamical Systems}, 36(1):509, 2016.

\end{thebibliography}

\end{document}